\documentclass[12pt]{amsart}

\newcommand {\supplus}{\mathop{{\supset}\llap{\raise
0.5pt\hbox{\normalfont\small+}\hskip 0.5pt}}}
\newcommand {\subplus}{\mathop{{\subset}\llap{\raise
0.5pt\hbox{\normalfont\small+}\hskip 0.5pt}}}
\newcommand {\Cee}    {{\mathbb  C}}

\newcommand {\Nee}    {{\mathbb  N}}

\newcommand {\Ree}    {{\mathbb  R}}

\newcommand {\Zee}    {{\mathbb  Z}}

\newcommand {\fA}     {{\mathfrak{A}}}

\newcommand {\fder}   {{\mathfrak{der}}}   %

\newcommand {\fg}     {{\mathfrak{g}}}    %
\newcommand {\fG}     {{\mathfrak{G}}}    %
\newcommand {\fgl}    {{\mathfrak{gl}}}  %
\newcommand {\fh}     {{\mathfrak{h}}}

\newcommand {\fo}     {{\mathfrak{o}}}
\newcommand {\fosp}   {{\mathfrak{osp}}}
\newcommand {\fp}    {{\mathfrak{p}}}   %

\newcommand {\fS}     {{\mathfrak{S}}}

\newcommand {\fsl}    {{\mathfrak{sl}}}

\newcommand {\fsp}    {{\mathfrak{sp}}}

\newcommand {\cal} {\mathcal}

\newcommand {\cF}     {{\cal F}}

\def \opname#1#2%
  {\expandafter\newcommand \csname #1\endcsname {{\mathop{#2}\nolimits}}}

\newcommand{\rmname}[1]
  {\expandafter\newcommand \csname #1\endcsname {{\operatorname{#1}}}}

\newcommand{\rmnameii}[2]
  {\expandafter\newcommand \csname #1\endcsname {{\operatorname{#2}}}}

\rmname{act}
\rmname{Ad}
\rmname{Add}
\rmname{ad}
\rmname{Alt}
\rmname{alt}
\rmname{Ann}
\rmname{antidiag}
\rmname{Ber}
\rmname{ber}
\rmname{Br}
\rmname{card}
\rmname{ch}
\rmname{Char}
\rmname{cem}
\rmname{cj}
\rmname{Cliff}
\rmname{cntr}
\rmname{codim}
\rmname{coind}
\rmname{const}
\rmname{col}
\rmname{cork}
\rmname{cpr}
\rmname{diag}
\rmnameii{Div}{div}
\rmname{Def}
\rmname{Der}
\rmname{Diff}
\rmname{Dim}
\rmname{End}
\rmname{Even}
\rmname{Ext}
\rmname{gr}
\rmname{Hom}
\rmname{HT}
\rmnameii{Ht}{ht}
\rmname{hwt}
\rmname{Id}
\rmname{id}
\rmname{ind}
\rmname{Ind}
\rmname{Inf}
\rmname{irr}
\rmname{Le}
\rmname{Lie}
\rmname{lwt}
\rmname{mult}
\rmname{Mat}
\rmname{Mor}
\rmname{nm}
\rmname{Ob}
\rmname{Odd}
\rmname{Osc}
\rmname{per}
\rmname{Pic}
\rmname{pr}
\rmname{pro}
\rmname{Prime}
\rmname{Proj}
\rmname{prt}
\rmname{pt}
\rmname{Q}
\rmname{qet}
\rmname{qtr}
\rmname{rd}
\rmname{rk}
\rmname{row}
\rmname{Res}
\rmname{salt}
\rmname{Sch}
\rmname{SBr}
\rmname{scalar}
\rmname{Ser}
\rmname{sign}
\rmname{Smbl}
\rmname{spin}
\rmname{ssym}
\rmname{str}
\rmname{st}
\rmname{sgn}
\rmname{sq}
\rmname{symm}
\rmname{supp}
\rmname{Supp}
\rmname{St}
\rmname{Spec}
\rmname{Spm}
\rmname{tr}
\rmname{vpt}
\rmname{weyl}
\rmname{Weyl}
\rmname{Witt}

\opname{vvol}  {{v\hspace{-0.1ex}o\hspace{-0.02ex}l\/}}
\opname{pnt}  {\text{\normalfont pt}}
\opname{Span} {{Span}}
\opname{slim} {\overline{\lim}}
\opname{Vol}  {{V\hspace{-0.55ex}o\hspace{-0.02ex}l\/}}
\opname{Par}  {{P\hspace{-0.3ex}a\hspace{-0.05ex}r\/}}

\newcommand {\eps} {\varepsilon}

\newcommand {\tto} {\longrightarrow}
\newcommand {\pder}[1] {{\frac{\partial}{\partial {#1}}}}

\newcommand {\bcdot}   {\mathbin{\hbox{\raise.4ex\hbox{\bf.}}}} % bold \cdot

\newcommand {\secno} {}
\newcommand {\ssecfont} {\normalfont\bf}

\newtheorem{Theorem}{\secno Theorem}

\newtheorem{Lemma}[Theorem]{\secno Lemma}
\newtheorem{Proposition}[Theorem]{\secno Proposition}
\newtheorem{Corollary}[Theorem]{\secno Corollary}

\newtheorem{Problem}[Theorem]{\secno Problem}

\newenvironment {th*}[1]
    {\gdef\thname{#1} \begin{thn}}%
    {\end{thn}}
\newtheorem{thn}[Theorem] {\thname}

\theoremstyle{definition}

\newenvironment {ex*}[1]
    {\gdef\thname{#1} \begin{exn}}%
    {\end{exn}}
\newtheorem{exn}[Theorem]{\thname}

\theoremstyle{remark}

\newtheorem{Remark}[Theorem]{\secno Remark}

\newenvironment {rem*}[1]
    {\gdef\thname{#1} \begin{remn}}%
    {\end{remn}}
\newtheorem{remn}[Theorem]{\thname}

\newcommand {\ssec}{\subsection*}

\newcommand {\ssbegin}[2]
  {\def \secno {\gdef \secno {}{\ssecfont #1. }}%
   \begin{#2}}
\begin{document}

\title[Orthogonal polynomials and complex size
matrices]{Orthogonal polynomials of discrete variable and Lie
algebras of complex size matrices}

\author{Dimitry Leites, Alexander Sergeev${}^1$}

\address{Correspondence: Department of Mathematics, University of Stockholm,
Roslagsv.  101, Kr\"aftriket hus 6, S-104 05, Stockholm, Sweden\\
mleites@math.su.se (${}^1$On leave of absence from Balakovo
University of Technique, Technology and Control, Balakovo, Saratov
region, Russia)}

\thanks{We are thankful to NFR for financial support, M.~Vasililev for
a stimulating question, P.~Grozman and, especially, V.~Gerdt for help
(sec.9.3).}



\begin{abstract} We give a uniform interpretation of the classical
continuous Chebyshev's and Hahn's orthogonal polynomials of
discrete variable in terms of Feigin's Lie algebra $\fgl(\lambda)$
for $\lambda\in\Cee$.  One can similarly interpret Chebyshev's and
Hahn's $q$-polynomials and introduce orthogonal polynomials
corresponding to Lie superalgebras.

We also describe the real forms of $\fgl(\lambda)$, quasi-finite
modules over $\fgl(\lambda)$, and conditions for unitarity of the
quasi-finite modules.  Analogs of tensors over $\fgl(\lambda)$ are
also introduced.
\end{abstract}

\maketitle

This is a transcript of the talk at MPI, Bonn in the memory of
Misha Saveliev (preprint MPI-1999-44 at www.mpim-bonn.mpg.de), see
also Theor.  Math.  Phys., v.  123, no.  2, 2000, 205--236
(Russian), 582--609 (English).  It consists of three parts: the
description of orthogonal polynomials proper (\S 1, 8, 9) and some
auxiliary results: the description of the generating function of
the trace and its analogs (\S 2) and a description of quasi-finite
modules over $\fgl(\lambda)$ for $\lambda\in\Cee$.

For a continuation, see \cite{S}.

\section*{\S 1. Introduction}

For recapitulations, see \cite{NU}, \cite{NSU}.

\ssec{1.0. Definitions} The equation of the form
$$
\sigma(x)y''+\tau(x) y'+\lambda y=0,\eqno{(1.0.1)}
$$
where $\sigma(x)$ is a polynomial of degree $2$, $\tau(x)$ is a
polynomial of degree $1$ and $\lambda$ is a constant, is called a
{\it hypergeometric type equation} and its solutions are called
functions of {\it hypergeometric type}.  By multiplying (0.1) by
an appropriate $\rho(x)$ one reduces (1.0.1) to the self-adjoint
form
$$
(\sigma(x)\rho(x)y')'+\lambda \rho(x)y=0,\text{ where
$(\sigma(x)\rho(x))'=\tau(x)\rho(x)$}.\eqno{(1.0.2)}
$$
Let $y_{m}$ and $y_{n}$ be solutions of (1.0.2) with distinct
eigenvalues $\lambda_{m}$ and $\lambda_{n}$, respectively.  If,
for some $a$ and $b$, not necessarily finite, $\rho(x)$ satisfies
the conditions
$$
\sigma(x)\rho(x)x^k|_{x=a, b}=0\text{ for } k=0, 1, \ldots ,
$$
then
$$
\int\limits_{a}^by_{m}(x)y_{n}(x)\rho(x)dx=0.\eqno{(1.0.3)}
$$
(Clearly, if $a$ and $b$ are finite, it suffices to require that
$\sigma(x)\rho(x)|_{x=a, b}=0$.)

Example: the {\it Jacobi polynomials} $P_{n}^{(\alpha, \beta)}(x)$
defined for $\alpha, \beta>-1$ as polynomial solutions of (0.2)
for $\sigma(x)=1-x^2$, $\rho=(1-x)^\alpha(1+x)^\beta$, $(a,
b)=(-1, 1)$ and $\tau=\beta-\alpha-(\alpha+\beta+2)x$.  For
$\alpha=\beta=0$, the Jacobi polynomials are called {\it
Chebyshev} polynomials.

The difference equation which approximates (1.0.1) on the uniform
lattice is, clearly, of the form
$$
\renewcommand{\arraystretch}{1.4}
\begin{array}{l}
\sigma(x)\frac{1}{h}\left[\frac{y(x+h)-y(x)}{h}-
\frac{y(x)-y(x-h)}{h}\right]
+\\
\frac{\tau(x)}{2}
\left[\frac{y(x+h)-y(x)}{h}+\frac{y(x)-y(x-h)}{h}\right]+\lambda
y=0.
\end{array}\eqno{(1.0.4)}
$$
Set $h=1$ and $\Delta f(x)=f(x+1)-f(x)$, $\nabla f(x)=f(x)-f(x-1)$.
Then the difference equation (1.0.4) takes the form
$$
\Delta(\sigma(x)\rho(x)\nabla y)+\lambda \rho(x)y=0,\text{ where
$\Delta(\sigma(x)\rho(x))=\tau(x)\rho(x)$}\eqno{(1.0.5)}
$$
while the orthogonality relations take the form
$$
\sum_{x_{i}=a}^{b-1}y_{m}(x_{i})y_{n}(x_{i})\rho(x_{i})=0\eqno{(1.0.6)}
$$
provided $\rho(x)$ satisfies the conditions
$$
\sigma(x)\rho(x)x^k|_{x=a, b}=0\text{ for } k=0, 1, \ldots
$$

{\bf Examples}: The {\it Hahn polynomial} $h_{n}^{(\alpha,
\beta)}(x, N)$ is a polynomial solution of (1.0.5) defined for
$\alpha, \beta>-1$ for $\sigma(x)=x(N+\alpha-x)$,
$\rho=\frac{\Gamma(N+\alpha-x)\Gamma(\beta+1+x)}{\Gamma(x+1)
\Gamma(N-x)}$, $(a, b)=(0, N)$ and
$\tau=(\beta+1)(N-1)-(\alpha+\beta+2)x$.  (As noted in \cite{R},
actually, Hahn's polynomials were known and studied by Chebyshev;
Hahn rediscovered them together with their $q$-analogs.)

The {\it Chebyshev polynomial} $t_{n}(x, N)$ is a polynomial solution
of (1.0.5) defined for $\alpha=\beta=0$, $\sigma(x)=x(N-x)$, $\rho=1$,
$(a, b)=(0, N)$ and $\tau=N-1-2x$.

A usual way to obtain $q$-analogs of these polynomials is to
consider nonuniform partitions of the segment $[a, b]$, see \cite{NSU}.
Our scheme applied to $U_{q}(\fsl(2))$ instead of  $U(\fsl(2))$
considered here gives another approach.

Until recently, the conventional inner products (1.0.3) and
(1.0.6) were supposed to be positive definite.  Though $N$, the
number of nodes of the uniform partition of segment $(a, b)$ is
integer, it is possible to replace it with any complex number in
expressions for $h_{n}^{(\alpha, \beta)}(x, N)$ and $t_{n}(x, N)$,
since the latter analytically depend on $N$.  It was observed (see
\cite{NSU}) that, for $N$ purely imaginary, there is a measure
leading to a positive definite inner product of type (1.0.3).  The
discovery of this fact reflects the luck and serendipity of the
researchers since no explanation of the phenomenon was seen.

It was also due to ingenious calculations that various identities
orthogonal polynomials satisfy were discovered; the completeness of
the list of such identities was never discussed.

\ssec{1.1.  Our result: a summary} We observed that the definition
and all identities the continuous Chebyshev and Hahn orthogonal
polynomials of discrete variable satisfy (and their versions on
inhomogeneous lattice, or, in modern terms, $q$-analogs) follow
from the properties of the Casimir elements for the Lie algebra
$\fgl(\lambda)$ of matrices of complex size (and its $q$-analog).
Our starting point was an attempt to explicitly calculate the
value of the quadratic Casimir operator on quasi-finite
$\fgl(\lambda)$-modules
--- a source of orthogonality relations and various identities the
said polynomials satisfy.  (For a more detailed than ours
description of quasi-finite modules of level 1, see \cite{Sh}.)

Let, first, $\lambda=n>0$ be an integer.  The invariant
functional, trace, on $\fG=\fgl(n)$ gives rise to the
non-degenerate bilinear form $(A, B)=\tr~AB$.  The homogeneous
components $\fG_{i}$, where $\fG_{i}$ is the space of matrices
with support on the $i$th over- (for $i>0$) and under-diagonals
(for $i<0$), are orthogonal with respect to this form and the form
is non-degenerate on $\fG_{-i}\oplus\fG_{i}$ and on $\fG_{0}$.
The reduction of the form to the canonical form on these spaces
leads to the classical Hahn and Chebyshev polynomials of discrete
variable, respectively.  By expressing the Casimir elements ---
the generators of the center of $U(\fG)$ --- in terms of these
polynomials in the standard basis one gets all the identities the
polynomials satisfy.  More exactly, we express elements of
$\fG_{0}$, as well as those of $\fG_{-i}\oplus\fG_{i}$, in terms
of polynomials in one variable, $H$. The details are given in the
main text; here is a gist of the idea.

In $\fsl(2)$, consider the basis
$Y=\begin{pmatrix}0&0\\1&0\end{pmatrix}$,
$H=\begin{pmatrix}1&0\\0&-1\end{pmatrix}$,
$X=\begin{pmatrix}0&1\\0&0\end{pmatrix}$ with the commutation
relations
$$
[X, Y]=H, \quad [H, Y]=-2Y, \quad [H, X]=2X.\eqno{(1.1.1)}
$$
Consider the {\it principal} embedding of $\fsl(2)$ into $\fG=\fgl(n)$
(this embedding corresponds to the $n$-dimensional irreducible
$\fsl(2)$-module):
$$
Y\mapsto \mathop{\sum}\limits_{i=1}^{n-1}E_{i+1, i},\;
H\mapsto \mathop{\sum}\limits_{i=1}^n(n-2i+1)E_{ii},\;
X\mapsto \mathop{\sum}\limits_{i=1}^{n-1}i(n-i)E_{i, i+1}.\eqno{(1.1.6)}
$$
Then $\fG$ is of the form $\fG=\mathop{\oplus}\limits_{|i|\leq
n-1}\fG_{i}$, where $\fG_{i}=\{X^i\fG_{0}\}$ and
$\fG_{-i}=\{\fG_{0}Y^i\}$ for $i>0$ and $\fG_{0}=\Cee[H]/(P_n(H))$ is
the Cartan subalgebra; here $P_n(H)=\mathop{\prod}\limits_{1\leq i\leq
n}(H-n+2i-1)$.

By applying this approach to $\fgl(\lambda)$, $\lambda\in \Cee$, we
clarify the known results on the co-called {\it continuous} versions
of the classical Hahn and Chebyshev polynomials of discrete variable,
cf.  \cite{NU}, \cite{NSU}, concerning purely imaginary values of
$N=\lambda$ and get some new ones.  So let us explain, first, what is
$\fgl(\lambda)$, the Lie algebra B.~Feigin introduced in \cite{F}.

The quadratic Casimir operator of $\fsl(2)$
$$
\Omega =2YX+\frac12 H^2+ H\eqno{(1.1.2)}
$$
lies in the center of $U(\fsl(2))$.  Let $I_{\lambda}$ be the
two-sided ideal in the associative algebra $U(\fsl(2))$ generated by
$\Omega-\frac12(\lambda^{2}-1)$.  It turns out that the associative
algebra $\tilde\fA_{\lambda}=U(\fsl(2))/I_{\lambda}$ is simple for
$\lambda\not\in \Zee\setminus\{0\}$, otherwise $\tilde\fA_{\lambda}$
contains an ideal such that the quotient is isomorphic to the matrix
algebra $\Mat(|\lambda|)$.  Set (\cite{Di1})
$$
\fA_{\lambda}=\left\{\begin{matrix}\tilde\fA_{\lambda}&\text{ if }
\lambda\not\in \Zee\setminus\{0\}\cr
\Mat(|\lambda|)&\text{ otherwise}.\end{matrix}\right.\eqno{(1.1.3)}
$$
As associative algebra with unit, $\fA_{\lambda}$ is generated by $X$,
$Y$ and $H$ subject to the relations (1.1.1) and
$$
XY=\frac14(\lambda^2-(H-1)^2) \eqno{(1.1.4)}
$$
and one more relation for integer values of $\lambda$:
$$
X^{|\lambda|}=0\text{  if }\lambda\in
\Zee\setminus\{0\}.\eqno{(1.1.5)}
$$
In what follows we set $\fgl(\lambda):=L(\fA_{\lambda})$, the Lie
algebra associated with the associative algebra $\fA_{\lambda}$ when we
replace the dot product with the bracket.

If in (1.3) we replace $U(\fsl(2))$ with its $q$-deformation,
$U_{q}(\fsl(2))$, we get $q$-versions of the classical polynomials,
identities, etc.  The details will be given elsewhere.

The classics considered orthogonal polynomials with respect to a
sign definite scalar product determined by a measure (or the
corresponding sum for the difference equations).  The existence of
such a measure for complex values of certain parameters looks like
a miracle in the traditional approaches, cf.  \cite{NU},
\cite{NSU}, \cite{VK}. Contrariwise, our approach makes it
manifest that there exists a measure for any
$\lambda\in\Zee\setminus\{0\}$ and for $\lambda$ purely imaginary,
see formula (9.2).

Moreover, we immediately see that (9.2) is sign definite for $\lambda$
real and such that $0<|\lambda|<1$.

We also see that since at generic point the scalar product is
INDEFINITE, it is natural to consider such scalar products as
well. Moreover, if we replace in the above construction $\fsl(2)$
with a Lie superalgebra $\fg$ we never have a sign definite scalar
product on its Cartan subalgebra.  Indefinite (but non-degenerate
and symmetric) scalar products were considered in the literature,
but the corresponding orthogonal polynomials satisfy differential
equations of order $>2$, cf.  \cite{LK}, \cite{MK}, \cite{D}. (We
are thankful to T.~Ya.~Azizov for these references.)  Such
equations describe, perhaps, some reality; this is plausible due
to the natural way they appear; M.~Vasiliev's ideas \cite{V} show
a possible scope of their applicability.  Still, modern physicists
prefer 2nd order equations, as describing processes more readily
at hand.  Though this attitude to solutions of higher order
equations might slacken with progress of science, we are glad that
polynomials orthogonal with respect to the restriction of $\tr$ on
$\fgl(\lambda)$ (and their superanalogs) satisfy a 2nd order
difference relation.

\ssec{1.2.  Several variables} In the known to us attempts to
generalize the classical continuous polynomials of discrete variable
to several variables one actually considers the representation of the
direct sum of several copies of $\fA_{\lambda}$, see (1.3) in the
tensor products of Verma modules over $\fsl(2)$.  This leads to the
sums of products of polynomials in one variable; more exactly, to
slightly more involved polynomials, but still, compositions of
polynomials in one variable.

Our approach leads to intrinsically much more involved polynomials,
namely, we suggest: in the above scheme, replace $\fsl(2)$ with ANY
Lie algebra or Lie superalgebra $\fg$.  This is possible provided
$\fg$ possesses (a) Verma modules with finitely many generators and
(b) several finite dimensional representations.  The theory seems to
be richest if $\fg$ is simple or ``close'' to simple (nontrivial
central extension, the derivation algebra, etc.).  Such
generalizations of $\fgl(\lambda)$ appear in spin $>2$ models and
Calogero-Sutherland model, see \cite{V}.

The results obtained in this way will be considered in separate
papers.  Finally, recently one of us (A.~S.) observed how to embrace
ALL classical polynomials of discrete variable, not just Chebyshev and
Hahn ones) in a scheme slightly generalizing the one described here.
This generalization and a superization of the above scheme (with
$\fsl(2)$ replaced with $\fosp(1|2)$ and $\fsl(1|2)$ was delivered at
the Special Invited Lecture at Advanced Study Institute, Special
Functions 2000: Current Perspective and Future Directions Arizona
State University, Tempe, Arizona, U.S.A., May 29 to June 9, 2000.

\section*{\S 2. The trace formula}

\ssbegin{2.1}{Lemma} Let $A$ be an associative algebra generated
by the Lie algebra $\fg$, conisdered as a subspace, i.e., $A$ is a
quotient of $U(\fg)$.  Then $[A, A]=[\fg, A]$.
\end{Lemma}

\begin{proof} It suffices to show that $[x_{1}\dots x_{n}, a]\in [\fg,
A]$ for any $a\in A$ and $x_{1}, \dots , x_{n}\in \fg$.  The
induction on $n$: for $n=1$, the statement is obvious.  Let $n>1$.
Let us make use of the identity $[ab, c]=[a, bc]+[b, ca]$ true in
any associative algebra, cf.  \cite{M}.  Then
$$
[x_{1}(x_{2}\dots x_{n}), a]=[x_{1}, x_{2}\dots x_{n}a]\pm
[x_{2}\dots x_{n}, ax_{1}].
$$
By inductive hypothesis, the second summand in the right hand side
lies in $[\fg, A]$.  \end{proof}

\ssbegin{2.1.1}{Corollary} Let $\fg$ be a simple Lie
algebra or $\fosp(1|2n)$.  Then $U(\fg)=Z(U(\fg))\oplus [U(\fg), U(\fg)]$.
\end{Corollary}

\begin{proof} Let $U(\fg)=Z(U(\fg))\oplus
\left(\mathop{\oplus}\limits_{\lambda\neq
0}n_{\lambda}L^{\lambda}\right )$ be the decomposition into
irreducible $\fg$-modules with respect to the adjoint representation.
Since the modules $L^{\lambda}$ are irreducible, $[\fg,
U(\fg)]=\left(\mathop{\oplus}\limits_{\lambda\neq
0}n_{\lambda}L^{\lambda}\right )$.  (This argument does not work for
Lie superalgebras distinct from $\fosp(1|2n)$ due to the lack of
complete reducibility.)
\end{proof}

Hereafter $\rho$ is a half-sum of positive roots.

\ssbegin{2.1.2}{Corollary} Let $\fh$ be a Cartan subalgebra in $\fg$,
which is either a simple finite dimensional Lie algebra; let
$\lambda\in\fh^{*}$, and $M^{\lambda-\rho}$ the Verma module with
highest weight $\lambda$ (NOT $\lambda-\rho$); let $J_{\lambda}$ be
the (left) ideal of $U(\fg)$ equal to the kernel of the representation
of $U(\fg)$ in $M^{\lambda-\rho}$.  Set $\fA_{\lambda}=
U(\fg)/J_{\lambda}$.  Then $\fA_{\lambda}= \Cee\oplus [\fA_{\lambda},
\fA_{\lambda}]$.
\end{Corollary}

\begin{proof} By Lemma 2.1 $[\fA_{\lambda}, \fA_{\lambda}]=[\fg,
\fA_{\lambda}]$; hence, as in Corollary 2.1.1, $\fA_{\lambda}=
Z(\fA_{\lambda})\oplus [\fA_{\lambda}, \fA_{\lambda}]$.  But
$Z(\fA_{\lambda})$ is a homomorphic image of $Z(U(\fg))$; hence, is
equal to $\Cee$.
\end{proof}

\ssec{2.2} If $\fg$ is a simple Lie algebra and $\lambda$ is a highest
weight of a finite dimensional module, then $P(\lambda)=\dim
L^\lambda$ is a polynomial in $\lambda=(\lambda_{1}, \dots ,
\lambda_{n})$.  Consider $P(\lambda)$ for any $\lambda\in\fh^{*}$ such
that $\mathop{\prod}\limits_{i\neq j}(\lambda_{i}-\lambda_{j})\neq 0$.
By Corollary 2.1.2, there is a unique functional $\tr$ on
$\fA_{\lambda}$ such that $\tr(ab)=\tr(ba)$ and $\tr(1)=P(\lambda)$.

\begin{Lemma} Let $u\in U(\fg)$, $\tilde u\in \fA_{\lambda}$ its
image.  Then $\tr(\tilde u)$ is also a polynomial in $\lambda$.
\end{Lemma}

\begin{proof} Let $\sharp: Z(U(\fg))\oplus [U(\fg), U(\fg)]\tto
Z(U(\fg))$ be the natural projection.  By definition of $\sharp$ we
have:
$$
\tr (\tilde u)=\chi_{\lambda}(u^\sharp)P(\lambda)=
\varphi(u^\sharp)(\lambda)P(\lambda),\eqno{(2.2)}
$$
where $\varphi: Z(U(\fg))\tto \Cee[\fh]$ is the Harish-Chandra
homomorphism, see \cite{Di}.
\end{proof}

\ssbegin{2.3}{Theorem} Let $\fg$ be a simple Lie algebra, $\fh$ its
Cartan subalgebra.  Let us identify the space $(\Cee[\fh])^*$ with the
algebra of formal power series on indeterminates $t=(t_{1}, \dots ,
t_{n})$ by the following formula:
$$
f\mapsto \mathop{\sum}\limits_{\nu=(\nu_{1},
\dots , \nu_{n})}f\left (\frac{h_{1}^{\nu_{1}}}{\nu_{1}!}
\ldots \frac{h_{n}^{\nu_{n}}}{\nu_{n}!}\right)t_{1}^{\nu_{1}}
\ldots  t_{n}^{\nu_{n}}, \eqno{(2.3.1)}
$$
where $h_{1}$, \dots , $h_{n}$ is a basis of $\fh$ and
$\nu_{i}\in\Zee_{+}$ for all $i$.  Let $e^{\lambda}(t)=
e^{\lambda_{1}t_{1}+\dots +\lambda_{n}t_{n}}$.  {\em (Recall that the
$\lambda_{i}$ are supposed to be distinct.)}

Then to the functional $\tr$ there corresponds the series
$$
\psi (\lambda, t)=\frac{\mathop{\oplus}\limits_{w\in W}\eps(w)
e^{w(\lambda)}} {\mathop{\oplus}\limits_{w\in W}\eps(w) e^{w(\rho)}}
(t).  \eqno{(2.3.2)}
$$
\end{Theorem}

\begin{proof} By Lemma 1.2 $\tr \left (\frac{h_{1}^{\nu_{1}}}{\nu_{1}!}
\ldots \frac{h_{n}^{\nu_{n}}}{\nu_{n}!}\right)$ is a polynomial in
$\lambda$.  If $\lambda\in P_{++}$ (here $P_{++}$ is the set of
highest weights of finite dimensional modules), then
$$
\hat \psi(\lambda, t)=\mathop{\sum}\limits_{\nu=(\nu_{1},
\dots , \nu_{n})}\tr \left (\frac{h_{1}^{\nu_{1}}}{\nu_{1}!}
\ldots \frac{h_{n}^{\nu_{n}}}{\nu_{n}!}\right)t_{1}^{\nu_{1}}
\ldots  t_{n}^{\nu_{n}}
$$
coincides with $(2.3.2)$ thanks to the Weyl character formula.
But $\psi(\lambda, t)=\mathop{\sum}\limits_{\nu}P_{\nu}t^{\nu}$,
where $P_{\nu}$ are some polynomials.  Indeed, since
$P_{\nu}(\lambda)=\tr\left (\frac{h_{1}^{\nu_{1}}}{\nu_{1}!}
\ldots \frac{h_{n}^{\nu_{n}}}{\nu_{n}!}\right)$ for $\lambda\in
P_{++}$, this is true for any $\lambda\in\fh^*$.  Hence, $\hat
\psi(\lambda, t)= \psi(\lambda, t)$.
\end{proof}

\section*{ \S 3. Highest weights of quasi-finite modules over
$\fgl(\lambda)$}

Let $L^i$ be the irreducible $\fsl(2)$-module with highest weight $i$.
The Lie algebra $\fgl(\lambda)$ considered as $\fsl(2)$-module with
respect to the adjoint action of the latter is of the form
$$
\fgl(\lambda)=\begin{cases} L^0\oplus L^2\oplus \dots
L^{2|\lambda|-2}\text{ for $\lambda\in \Zee\setminus\{0\}$}\cr
L^0\oplus L^2\oplus \dots \text{ otherwise.
}\end{cases}\eqno{(3.0.1)}
$$

The Lie algebra $\fgl(\lambda)$ possesses a $\Zee$-grading
$\fgl(\lambda)=\mathop{\oplus}\limits_{i\in\Zee}\fgl(\lambda)_{i}$,
where
$$
\fgl(\lambda)_{i}=\{z\in\fgl(\lambda): [H, z]=2iz\}\eqno{(3.0.2)}
$$
and the increasing filtration:
$\fgl(\lambda)_{(i)}=\mathop{\oplus}\limits_{k\leq
i}\fgl(\lambda)_{k}$.  Therefore, we can speak about the {\it
triangular decomposition} $\fg=\fg_{-}\oplus\fg_{0}\oplus\fg_{+}$ of
$\fg=\fgl(\lambda)$, where $\fg_{\pm}=\mathop{\oplus}\limits_{i\geq
0}\fgl(\lambda)_{\pm i}$, as well as about {\it parabolic subalgebras}, {\it highest
weights}, etc., see \cite{KR}.

Following \cite{KR} we say that the $\fg$-module $V$ is {\it
quasi-finite} if $V=\mathop{\oplus}\limits_{j\in\Zee}V_{j}$ and
$\dim V_{j}<\infty$ for every $j$.  (We only consider {\it graded}
$\fg$-modules, i.e., such that $\fg_{i}V_{j}\subset V_{i+j}$.)

The well-known fact that $\fgl(n)$ has a nontrivial one-dimensional
representation ($A\mapsto \tr~A$) has its counterpart for
$\fgl(\lambda)$.

\ssec{3.1} A subalgebra $\fp\subset \fg$ is called a {\it
parabolic} one if it contains $\fg_0\oplus \fg_+$ as a proper
subalgebra. For example, for any $P\in\Cee[H]$, set
$\fp_{-1}(P)=P\Cee[H]Y$ and let $\fp(P)$ be generated by
$\fp_{-1}(P)$ and $\fg_0\oplus \fg_+$. This is the {\it minimal}
parabolic subalgebra corresponding to $P$.

\begin{Lemma} The minimal parabolic subalgebra $\fp(P)$ is of
the form $\fp(P)=\mathop{\oplus}\limits_{k\in\Zee}\fp(P)_k$, where
$\fp(P)_k=\fg_k$ for $k\geq 0$ and $\fp(P)_{-k}=I_k\cdot Y^ k$, where
$I_k$ is the ideal of $\Cee[H]$ generated by $P(H)P(H+2)\dots
P(H+2k-2)$.
\end{Lemma}

\begin{proof} See \cite{KR}, mutatis mutandis.\end{proof}

\ssec{3.2} Let $\lambda\in\fg^{*}_{0}$, $\fg_{+}v_{\lambda}=0$ and
$M^\lambda=\ind^{\fg}_{\fg_{0}\oplus \fg_{+}}(\Cee v_{\lambda})$ the
corresponding Verma module; let $L^{\lambda}$ be the irreducible
$\fg$-module with highest weight $\lambda$.

\begin{Lemma} {\em (\cite{KR})} The following
conditions are equivalent:

{\em i)} $M^ {\lambda}$ contains a vacuum vector that lies in $M^
{\lambda}_{(-1)}$, the space of filtration $-1$;

{\em ii)} The module $L^ {\lambda}$ is quasi-finite;

{\em iii)} $L^ {\lambda}$ is a quotient of the generalized Verma
module $M^ {\lambda, P}=\ind^{\fg}_{\fp(P)}(\Cee v_{\lambda})$ for
some $P\in\Cee[H]$.
\end{Lemma}

\ssec{3.3} Shoikhet calls $\deg~P$ the {\it level} of the module $M^
{\lambda, P}$.  In \cite{Sh} he described modules of level 1 in more
details than we give here.

We have already established that
$\fg^*_0\simeq(\Cee[H])^*\simeq\Cee[[t]]$, the algebra
isomorphism being given by the map $F$:
$$
F: (\Cee[H])^*\tto \Cee[[t]], \quad\theta\mapsto
F_\theta(t)=\mathop{\sum}\limits_{k=0}^\infty\frac{\theta(H^k)}{k!}t^k.
$$

The following theorem describes the set of formal power series
corresponding to the highest weights of quasi-finite modules.
(Recall that a {\it quasi-polynomial} is an expression of the form
$\sum R_i(t)e^{\alpha_{i}t}$, where the $R_i(t)$ are polynomials.)

\begin{Theorem} The formal power series corresponding to the highest
weight of any quasi-finite module over $\fgl(\lambda)$ is of the
form $\displaystyle\frac{R(t)}{1-e^{-2t}}$, where $R(t)$ is a
quasi-polynomial such that $R(0)=0$.  \end{Theorem}

\begin{proof} First, observe the following easy to verify statements:
(a) if $\theta\in(\Cee[H])^*$ and
$F_\theta(t)=\mathop{\sum}\limits_{k=0}^\infty\frac{\theta(H^k)}{k!}t^k$
the corresponding series, then
$$
\mathop{\sum}\limits_{k=0}^\infty\frac{\theta((H+a)^k)}{k!}t^k=
e^{at}F_\theta(t);
$$
and (b) for any  $R(H)\in\Cee[H]$ we have
$$
\mathop{\sum}\limits_{k=0}^\infty\frac{\theta(R(H)H^k)}{k!}t^k=
R(\frac{d}{dt})F_\theta(t)\text{ and }
\mathop{\sum}\limits_{k=0}^\infty\frac{\theta(R(H+a)(H+a)^k)}{k!}t^k=
R(\frac{d}{dt})(e^{at}F_\theta(t)).
$$
If $\theta$ is the highest weight of a quasi-finite module over
$\fgl(\lambda)$, then by Lemma 3.2 there exists a polynomial
$P(H)\in\Cee[H]$ such that $\theta$ can be extended to a
one-dimensional representation of the minimal parabolic subalgebra
corresponding to $P$.  It is not difficult to verify that $[\fp,
\fp]\cap\fg_0=[\fg_{1}, \fp_{-1}]$.  Hence, $\theta([\fg_{1},
\fp_{-1}])=0$.

Denote:
$$
T(H)=XY=\frac14(\lambda^2-(H+1)^2).\eqno{(3.2)}
$$
Then the condition $\theta([\fg_{1}, \fp_{-1}])=0$ can be
expressed as $\theta([X, P(H)H^kY])=0$, or, equivalently, as
$\theta(XP(H)H^kY-P(H)H^kYX])=0$, or, as
$$
\theta(T(H-2)P(H-2)(H-2)^k - P(H)H^kT(H))=0.
$$
Therefore,
$$
\renewcommand{\arraystretch}{1.4}
\begin{array}{l}
    \mathop{\sum}\limits_{k=0}^\infty\displaystyle
    \frac{\theta(T(H-2)P(H-2)(H-2)^k)}{k!}t^k-
    \mathop{\sum}\limits_{k=0}^\infty\displaystyle
    \frac{\theta(T(H)P(H)(H)^k)}{k!}t^k=\\
T(\frac{d}{dt})P(\frac{d}{dt})(e^{-2t}F_\theta(t))-
T(\frac{d}{dt})P(\frac{d}{dt})(F_\theta(t))=\\
T(\frac{d}{dt})P(\frac{d}{dt})\left((e^{-2t}-1)F_\theta(t)\right)=0.
\end{array}
$$
Thus, the function $(e^{-2t}-1)F_\theta(t)$ is a solution of an
ordinary differential equation with constant coefficients, hence,
is a quasi-polynomial.  Obviously, $R(0)=0$.

Conversely, if $R(t)$ is a quasi-polynomial, $R(0)=0$ and
$P(\frac{d}{dt})R(t)=0$, then
$T(\frac{d}{dt})P(\frac{d}{dt})R(t)=0$ and, by setting
$F_\theta(t)=\displaystyle\frac{R(t)}{e^{-2t}-1}$ (which is well
defined since $R(0)=0$) we satisfy the condition $\theta([\fg_{1},
\fp_{-1}])=0$ for $\fp_{-1}=P(H)\fg_0Y$.  Hence, $\theta$ is the
highest weight of a quasi-finite module.  \end{proof}

\ssec{3.4.  The trace on $\fgl(0)$} If $P=1$, then the parabolic
subalgebra of $\fg=\fgl(0)$ coincides with the whole algebra and,
therefore, $[\fg, \fg]\cap\fg_0=[\fg_1, \fg_{-1}]\cap\fg_0\not\ni 1$.
This proves that $\fg\neq [\fg, \fg]$ and, therefore, proves also
existence of an invariant functional $\theta$ on $\fg=\fgl(\lambda)$.
Therefore, $R=(1-e^{-2t})F_\theta(t)$ satisfies the equation
$T(\frac{d}{dt})R(t)=0$; explicitly:
$$
(\lambda^ 2-(\frac{d}{dt}+1)^ 2)R(t)=0.
$$
If $\lambda\neq 0$, then solutions are of the form
$R(t)=c_1e^{(\lambda-1)t}+c_2e^{-(\lambda+1)t}$ and the
initial condition $R(0)=0$ implies that
$$
F_\theta(t)=
c\frac{e^{(\lambda-1)t}-e^{-(\lambda+1)t}}{1-e^{-2t}}.
$$
Let us identify the scalars with the corresponding scalar matrices and norm the
functional $\theta$ (trace) naturally, i.e., by setting
$\theta(1)=\lambda$. This fixes $c$, namely, $c=1$.

Clearly, if $\lambda=0$, then the characteristic equation has
multiple roots and
$$
F_\theta(t)=c\frac{te^{-t}}{1-e^{-2t}}.
$$
Thus, $\fsl(0)$ is an infinite dimensional Lie algebra
with an $\infty$-dimensional identity module.

\section*{\S 4. The Lie algebras \protect $\fgl(\lambda)$ and
$\fgl^-(\infty)$, $\fgl^+(\infty)$, $\fgl(\infty)$}

Let $V$ be a vector space with a fixed basis $v_i$ for $i\in\Zee$, let
$V^+$ be its subspace spanned by $v_i$ for $i\geq 0$ and $V^-$ its
subspace spanned by $v_i$ for $i<0$.  Let $\fgl^-(\infty)$,
$\fgl^+(\infty)$, and $\fgl(\infty)$ be the Lie algebras of linear
transformations of $V^-$, $V^+$ and $V$, respectively, whose matrices
in the fixed bases are supported on a finite number (depending on the
matrix) of diagonals parallel to the main one.

On  $\fgl(\infty)$, there is a (unique) nontrivial 2-cocycle
$$
c(A, B)=\tr~([J, A]B),\;\text{ where } J=\mathop{\sum}\limits_{i\leq 0}E_{ii}
-\mathop{\sum}\limits_{i> 0}E_{ii}.\eqno{(4.0)}
$$
Denote the corresponding central extension by
$\widehat{\fgl}(\infty)$; the bracket in $\widehat{\fgl}(\infty)$ is
of the form
$$
[A, B]=AB-BA+\tr~([J, A]B)\cdot z,\text{ where $z$ is the new central
element}.
$$
In particular, for matrix units, we have
$$
[E_{ij}, E_{kl}]=\delta_{jk}E_{il}-\delta_{il}E_{kj}+
\delta_{il}\delta_{jk}(\kappa(i)-\kappa(j))\cdot z,
$$
where $\kappa(i)=1$ if $i\leq 0$ and $0$ otherwise.

\ssec{4.0.  Harish-Chandra modules over $\fsl(2)$}
Harish-Chandra modules $M^{\lambda, s}=\Span(v_{i}: i\in\Zee)$ for any
$\lambda, s\in\Cee$ are defined over $\fsl(2)$ as the span of the
$v_i$ that satisfy the following relations
$$
\renewcommand{\arraystretch}{1.4}
\begin{array}{l}
Hv_{i}=(s+2i)v_{i},\\
Xv_{i}=\frac12\sqrt{(\lambda-s-2i-1)(\lambda+s+2i+1)}v_{i+1},\\
Yv_{i}=\frac12\sqrt{(\lambda-s-2i+1)(\lambda+s+2i-1)}v_{i-1}.
\end{array}\eqno{(4.0.1)}
$$
It is not difficult to observe that the quadratic Casimir operator
$\Omega$ acts on $M^{\lambda, s}$ as a scalar operator of
multiplication by $\frac 12(\lambda^{2}-1)$ and $M^{\lambda, s}\simeq
M^{\lambda, s'}$ if $s-s'\in 2\Zee$.

Suppose $\lambda\in \Zee\setminus \{0\}$.  Then by \cite{Di} if
$\lambda-s\not\in 2\Zee+1$, then $M^{\lambda, s}$ is irreducible,
whereas if $\lambda-s\in 2\Zee+1$, then $M^{\lambda, s}=M^{\lambda,
s}_{+}\oplus M^{\lambda, s}_-$, where $M^{\lambda, s}_+$ is an
irreducible $\fsl(2)$-module with the lowest weight $\lambda+1$ and
$M^{\lambda, s}_-$ is an irreducible $\fsl(2)$-module with the highest
weight $\lambda-1$ and these module exhaust all irreducible
$\fsl(2)$-modules with a diagonal $H$-action on which $\Omega$ acts as
a scalar operator of multiplication by $\frac12(\lambda^{2}-1)$.

There is a more transparent realization of Harish-Chandra modules,
namely, set
$$
X^+=x\pder{y}-2(\lambda-1)\frac{x}{y},
\;X^-=y\pder{x};\;\text{(hence, $H=x\pder{x}-y\pder{y}-(\lambda-1)$)}
$$
and set $w_i=x^{a+i}y^{b-i}$.  Then, $\lambda=a+b+1$, $s=a-b$ and
$w_i=v_{i}$ up to a scalar factor.

\ssec{4.1} The module $M^{\lambda, s}$ over $\fgl(\lambda)$ determines
an embedding $\phi: \fgl(\lambda)\tto \fgl(\infty)$.  Since
$H^2(\fgl(\lambda))=0$ (see \cite{F}), there exists a lift of the
embedding $\phi$ to a homomorphism $\hat\phi: \fgl(\lambda)\tto
\widehat{\fgl}(\infty)$ which is of the form
$\hat\phi(u)=\phi(u)+\theta(u)\cdot z$, where $\theta$ is a linear
functional on $\fgl(\lambda)$ such that $\theta([u, v])=c(\phi(u),
\phi(v))$ for the cocycle on $\fgl(\infty)$ which makes it into
$\widehat{\fgl}(\infty)$.  Observe that $\theta$ is determined
uniquely up to a functional proportional to the trace on
$\fgl(\lambda)$.

The following theorem describes $\hat\phi$ more explicitly.

\begin{Theorem} Let $\fg=\fgl(\lambda)$.  Then
$\hat\phi(\fg_i)=\phi(\fg_i)$ for $i\neq 0$ and
$$
\hat\phi(e^{tH})=\phi(e^{tH})-\frac{e^{st}-
e^{-(\lambda+1)t}}{1-e^{-2t}}z.\eqno{(4.1.1)}
$$
\end{Theorem}

\begin{proof}
$$
\renewcommand{\arraystretch}{1.4}
\begin{array}{l}
\phi(X)=\sum\alpha_iE_{i+1, i},\text{ where
}\alpha_i=\frac12\sqrt{(\lambda-s-2i-1)(\lambda+s+2i+1)},\\
\phi(H)=\sum\gamma_iE_{i, i},\text{ where
}\gamma_i=s+2i,\\
\phi(Y)=\sum\beta_iE_{i-1, i},\text{ where
}\beta_i=\frac12\sqrt{(\lambda-s-2i+1)(\lambda+s+2i-1)}.\\
\end{array}\eqno{(4.1.2)}
$$
Therefore,
$$
[J, \phi(X)]=\sum(\kappa(i)-\kappa(i+1))\alpha_iE_{i+1, i}
$$
and $\phi(H^k)\phi(Y)=\sum\gamma_i^k\beta_{i+1}E_{i, i+1}$.
Hence,
$$
[J, \phi(X)]\phi(H^k)\phi(Y)=\sum\gamma_i^k\beta_{i+1}
\alpha_i(\kappa(i)-\kappa(i+1))E_{i+1, i+1}
$$
and
$$
\renewcommand{\arraystretch}{1.4}
\begin{array}{l}
c(\phi(X), \phi(H^kY))=\tr([J, \phi(X)]\phi(H^kY))=
\sum\gamma_i^k\beta_{i+1}
\alpha_i(\kappa(i)-\kappa(i+1))=\\
\gamma_0^k\beta_{1} \alpha_0= \frac14s^k(\lambda^2-(s+1)^2)=T(s)s^k.
\end{array}
$$
Therefore,
$$
\renewcommand{\arraystretch}{1.4}
\begin{array}{l}
T(\frac{d}{dt})(e^{st})=T(s)e^{st}=\mathop{\sum}\limits_{k\geq
0}\displaystyle\frac{T(s)s^k}{k!}t^k =\\
\mathop{\sum}\limits_{k\geq 0}\displaystyle\frac{c(\phi(X), \phi(H^kY))}{k!}t^k =\\
\mathop{\sum}\limits_{k\geq 0}\displaystyle\frac{\theta([ X, H^kY])}{k!}t^k =
\mathop{\sum}\limits_{k\geq 0}\displaystyle\frac{\theta(XH^kY-H^kYX)}{k!}t^k=\\
\mathop{\sum}\limits_{k\geq 0}\displaystyle\frac{\theta(T(H-2)(H-2)^k-H^kT(H))}{k!}t^k=\\
T(\frac{d}{dt})\left((e^{-2t}-1)\theta(e^{tH})\right ).
\end{array}
$$
Therefore,
$T(\frac{d}{dt})\left((e^{-2t}-1)\theta(e^{tH})-e^{st})\right )=0$
and
$$
(e^{-2t}-1)\theta(e^{tH})=\begin{cases} e^{st}+
c_1e^{-(\lambda+1)t}+c_2e^{(\lambda-1)t}&\text{ if $\lambda\neq
0$}\cr e^{st}+c_1e^{-t}+c_2te^{-t}&\text{ if
$\lambda=0$.}\end{cases}
$$
In both cases, subtracting the summand proportional to the
trace we get
$$
(e^{-2t}-1)\theta(e^{tH})=e^{st}-e^{-(\lambda+1)t}.
$$
\end{proof}

\ssec{4.2.  $\fgl(\infty)$ over the algebra of truncated polynomials}
Let $R_m=\Cee[\eps]$, where $\eps^{m+1}=0$.  By extending the module
$M^{\lambda, s}$ to a module $M^{\lambda, s}_{R_{m}}$ over $R_m$,
i.e., by assuming that $s\in R_m$, we may define the action of
$\fsl(2)$ in $M^{\lambda, s}_{R_{m}}$  by the {\it same} formulas as
in $M^{\lambda, s}$ .  If the $v_i$ is the initial
basis, then for the fixed basis of the extended module $M^{\lambda,
s}_{R_{m}}$ we take $\eps^jv_i$ for all $i\in\Zee$ and $j=1$, \dots ,
$m$.

Observe that in this basis of $M^{\lambda, s}_{R_{m}}$ the action of
$H$ is not a diagonal one.

The algebra $\fgl(\infty; R_m)$ and its subalgebras
$\fgl^{\pm} (\infty; R_m)$  are naturally defined and the central
extension $\widehat{\fgl}(\infty; R_m)$ is determined by the
same formula.

By considering $s+\eps$ instead of $s\in\Cee$, we get the
following corollary of Theorem 4.1:

\begin{Corollary} The homomorphism $\hat\phi:\fgl(\lambda)\tto
\widehat{\fgl}(\infty; R_m)$ induced by the embedding
$\phi:\fgl(\lambda)\tto \fgl(\infty; R_m)$ satisfies
$\hat\phi(\fgl(\lambda)_i)=\phi(\fgl(\lambda)_i)$ for $i\neq 0$ and
$$
\hat\phi(e^{tH})=\phi(e^{tH})-
\left(\frac{e^{st}-e^{-(\lambda+1)t}}{1-e^{-2t}}+
\frac{e^{st}}{1-e^{-2t}}
\mathop{\sum}\limits_{j=1}^m\frac{\eps^{j}t^{j}}{j!}\right)z.\eqno{(4.2)}
$$
\end{Corollary}

\ssbegin{4.3}{Remark} In order to describe quasi-finite modules
over $\fgl(\lambda)$ we are ``forced" to consider the central
extension $\widehat{\fgl}(\infty)$ because
$\widehat{\fgl}(\infty)$ has more quasi-finite representations
than $\fgl(\infty)$ does and it turns out to be insufficient to
consider $\fgl(\infty)$-modules only.
\end{Remark}

\ssec{4.4.  The highest weights of quasi-finite modules} The
following statements describe the restrictions onto the highest
weight of an irreducible module necessary and sufficient for the
module to be quasi-finite.

Denote by $\fgl_f(\infty)$ and $\fgl_f^{\pm}(\infty)$ the subalgebras
of $\fgl(\infty)$ and $\fgl^{\pm}(\infty)$, respectively, consisting
of matrices with finite support.  Let $\fh$ denote the Cartan
subalgebra of any of these Lie algebras; consider it spanned by the
diagonal elements $E_{ii}$.  Observe that the restriction of the
embedding $\hat\phi: \fgl_f(\infty)\tto\widehat{\fgl}(\infty)$ onto
$\fh$ is of the form $E_{ii}\mapsto \phi(E_{ii})+\kappa(i)z$, where
$\phi$ is the natural embedding of $\fgl_f(\infty)$ into
$\fgl(\infty)$.  If $\theta$ is a weight of a (highest weight) module
over one of the above Lie algebras, both $\fgl(\infty)$ and its
subalgebras and $\fgl_f(\infty)$ and its subalgebras, we may,
considering $\fh$ in the above way, determine the coordinates
$\theta_i=\theta(E_{ii})$ of $\theta$.

\ssbegin{4.4.1}{Proposition} The $\fgl_f(\infty)$-module with
highest weight $\Lambda$ is quasi-finite if and only if there are
only finitely many distinct coordinates $\lambda_i$, where
$\lambda_i=\Lambda(E_{ii})$.
\end{Proposition}

\begin{proof} If $V$ is quasi-finite and generated by the
highest weight vector $v$. Since $\dim V_{-1}<\infty$, it
follows that  $V_{-1}=\Span(E_{i+1, i}v)_{i\in I}$ for
a finite set $I$. Let $k\not\in I$ and  $k+1\not\in I$. Then
$E_{k+1, k}v=\mathop{\sum}\limits_{i\in I}\alpha_i E_{i+1, i}v$, so
$E_{k, k+1}(E_{k+1, k}v)=0$ or, equivalently, $(E_{k, k}-E_{k+1,
k+1})v=0$, i.e.,  $\lambda_{k}=\lambda_{k+1}$.

The opposite implication is obvious. \end{proof}

\ssbegin{4.4.2}{Proposition} The $\fgl(\infty)$-module $V$ with
highest weight $\Lambda$, where $\lambda_i=\lambda(E_{ii})$, is
quasi-finite if and only if there are only finitely many nonzero
coordinates $\lambda_i$, where $\lambda_i=\Lambda(E_{ii})$.
\end{Proposition}

\begin{proof} By considering $V$ as $\fgl_f(\infty)$-module we deduce
that there are only finitely many distinct coordinates $\lambda_i$, so
we can set $\lambda_k=\lambda_+$ and $\lambda_{-k}=\lambda_-$ for
sufficiently large $k$.  Hence, $E_{k, k+1}(E_{k+1, k}v)=0$ and,
thanks to irreducibility of $V$, $E_{k+1, k}v=0$ for all $k$ but
finitely many.

Let $A=\mathop{\sum}\limits_{i\in\Zee} E_{i, i+1}$ and
$B_k=\mathop{\sum}\limits_{i\geq k} E_{i+1, i}$.  Then
$$
[A, B_k]=\mathop{\sum}\limits_{i\geq k}E_{i,
i}-\mathop{\sum}\limits_{i\geq k}E_{i+1, i+1}= E_{k, k}.
$$
Therefore,
$$
\renewcommand{\arraystretch}{1.4}
\begin{array}{l}
[A, B_k]v=E_{k, k}v=\lambda_kv=AB_kv-B_kAv=AB_kv=\\
A(\mathop{\sum}\limits_{i\geq k, \, i\in I}E_{i+1,
i}v)=\mathop{\sum}\limits_{i\geq k, \, i\in I}[A, E_{i+1,
i}]v=\mathop{\sum}\limits_{i\geq k, \, i\in I}(E_{i, i}-E_{i+1,
i+1})v=\\
\left(\mathop{\sum}\limits_{i\geq k, \, i\in
I}(\lambda_{i}-\lambda_{i+1})\right)v=(\lambda_k-\lambda_j)v,
\end{array}
$$
where $j$ is the largest index in the finite set $I$.
Therefore, $\lambda_k=\lambda_k-\lambda_j$ and
$\lambda_j=0$.

Similarly, if $C_k=\mathop{\sum}\limits_{i\leq k}E_{i+1, i}$, then
$$
[A, C_k]=\mathop{\sum}\limits_{i\leq k}E_{i,
i}-\mathop{\sum}\limits_{i\leq k}E_{i+1, i+1}=E_{k+1, k+1}.
$$
Hence,
$$
\renewcommand{\arraystretch}{1.4}
\begin{array}{l}
[A, C_k]v=\lambda_{k+1}v=AC_kv=A(\mathop{\sum}\limits_{i\leq k}E_{i+1,
i})v= \mathop{\sum}\limits_{i\leq k}[A, E_{i+1, i}]v\\
\mathop{\sum}\limits_{i\leq k}(E_{i, i}-E_{i+1, i+1})v=\mathop{\sum}\limits_{i\leq
k}(\lambda_{i}-\lambda_{i+1})v=(\lambda_{j}-\lambda_{k+1})v.
\end{array}
$$
So, $\lambda_{j}=0$ for sufficiently large $j$. The converse
statement follows from Proposition 4.4.1. \end{proof}

\ssbegin{4.4.3}{Proposition} The $\widehat{\fgl}(\infty)$-module
$V$ with highest weight $(\Lambda, c)$, where $c$ is the value of
the central charge on $z$, is quasi-finite if and only if there
are only finitely many nonzero coordinates $\lambda_i$, where
$\lambda_i=\Lambda(E_{ii})$.
\end{Proposition}

\begin{proof} As in the proof of Proposition 4.4.2 we show that there
are only finitely many distinct coordinates $\lambda_i$, and
$E_{k+1, k}v\neq 0$ for finitely many values of $k$.

Let $A=\mathop{\sum}\limits_{i\in\Zee} E_{i, i+1}$ and
$J=\mathop{\sum}\limits_{i\leq 0} E_{i, i}$.  Then $[A, J]=E_{01}$,
so $c(A, B)=\tr~(E_{01}B)$ for any $B\in\fgl(\infty)$, so $c(A,
B_k)=\kappa(k)$.  We have
$$
\renewcommand{\arraystretch}{1.4}
\begin{array}{l}
[A, B_k]v=E_{k, k}v=(E_{k, k}+\kappa(k)z)v=AB_kv=\\
A(\mathop{\sum}\limits_{i\geq k, \, i\in I}E_{i+1,
i}v)=\mathop{\sum}\limits_{i\geq k, \, i\in I}[A, E_{i+1,
i}]v=\mathop{\sum}\limits_{i\geq k, \, i\in I}(E_{i, i}-E_{i+1,
i+1}+\delta_{i0}z)v=\\
\kappa(k)z+\left(\mathop{\sum}\limits_{i\geq k, \, i\in
I}(\lambda_{i}-\lambda_{i+1})\right)v=\kappa(k)cv+
(\lambda_k-\lambda_+)v=(\kappa(k)c+
(\lambda_k)v,
\end{array}
$$
so $\lambda_+=0$.

For $C_k=\mathop{\sum}\limits_{i\leq k}E_{i+1, i}$, we have $c(A,
C_k)=\kappa(-k)$, so
$$
\renewcommand{\arraystretch}{1.4}
\begin{array}{l}
[A, C_k]v=(-E_{k+1, k+1}-\kappa(-k))v=\\
(-\lambda_{k+1}-\kappa(-k)c)v=AC_kv=A\mathop{\sum}\limits_{i\leq
k}E_{i+1, i}v=\mathop{\sum}\limits_{i\leq k}[A, E_{i+1, i}]v=\\
\mathop{\sum}\limits_{i\leq
k}(E_{i,i}-E_{i+1,i+1}+\delta_{i0}z)v=
(\kappa(-k)c+\lambda_{-}-\lambda_{k+1})v.
\end{array}
$$
So, $\lambda_{-}=0$. \end{proof}

\ssbegin{4.4.4}{Theorem} Let $\fg$ be one of the Lie algebras
$\fgl(\infty; R_m)$ or $\fgl^{\pm}(\infty; R_m)$ or their hatted
versions.  The $\fg$-module $V$ with highest weight $\theta$ given
by its coordinates $\theta_{ij}=(\eps^jE_{ii})$ is quasi-finite if
and only if there are only finitely many nonzero coordinates
$\theta_{ij}$.
\end{Theorem}

Proof is similar to that of 4.4.1--4.4.3.

Let $V$ be one of the modules $M^{\lambda, s}_{R_{m}}$ or, if
$M^{\lambda, s}_{R_{m}}$ is reducible, one of its irreducible
submodules.  We obtain a homomorphism of $\fgl(\lambda)$ into
$\fg=\fgl(\infty; R_m)$ or $\fgl^{\pm}(\infty; R_m)$.  Let $\theta$ be
a linear functional satisfying the conditions of Theorem 4.4.4.  We
may consider $\theta$ be a linear functional on the Cartan subalgebra
of $\fgl(\lambda)$.  Let us calculate the corresponding generating
function.

\ssbegin{4.5}{Theorem} The generating function $F_{\theta}(t)$
for the functional $\theta$ is as follows:

{\em i)} For $M^{\lambda, s}_{R_{m}}$:
$$
F_{\theta}(t)=\frac{\mathop{\sum}\limits_{i\in\Zee}e^{(s-2i)t}
\mathop{\sum}\limits_{j=0}^m
\displaystyle\frac{(\theta_{ij}-\theta_{i-1,
j})t^j}{j!}}{1-e^{-2t}}-
\left(\frac{e^{st}-e^{-(\lambda+1)t}}{1-e^{-2t}}c+
\frac{e^{st}}{1-e^{-2t}}
\mathop{\sum}\limits_{j=1}^m\frac{c_j}{j!}\right)
\eqno{(4.5.1)}
$$

{\em ii)} For the module with highest weight $\lambda-1$:
$$
F_{\theta}(t)=\frac{\mathop{\sum}\limits_{i\in\Zee}e^{(\lambda-2i-1)t}
\mathop{\sum}\limits_{j=0}^m
\displaystyle\frac{(\theta_{ij}-\theta_{i-1, j})t^j}{j!}}{1-e^{-2t}}
\eqno{(4.5.2)}
$$

{\em iii)} For the module with highest weight $-\lambda-1$:
$$
F_{\theta}(t)=\frac{\mathop{\sum}\limits_{i\in\Zee}e^{(-\lambda-2i-1)t}
\mathop{\sum}\limits_{j=0}^m
\displaystyle\frac{(\theta_{ij}-\theta_{i-1, j})t^j}{j!}}{1-e^{-2t}}
\eqno{(4.5.3)}
$$

{\em iv)} For the module with lowest weight $\lambda+1$:
$$
F_{\theta}(t)=\frac{\mathop{\sum}\limits_{i\in\Zee}e^{(\lambda+2i+1)t}
\mathop{\sum}\limits_{j=0}^m
\displaystyle\frac{(\theta_{ij}-\theta_{i-1, j})t^j}{j!}}{1-e^{-2t}}
\eqno{(4.5.4)}
$$

{\em v)} For the module with lowest weight $1-\lambda$:
$$
F_{\theta}(t)=\frac{\mathop{\sum}\limits_{i\in\Zee}e^{(1-\lambda+2i)t}
\mathop{\sum}\limits_{j=0}^m
\displaystyle\frac{(\theta_{ij}-\theta_{i-1, j})t^j}{j!}}{1-e^{-2t}}
\eqno{(4.5.5)}
$$
\end{Theorem}

\begin{proof} We will only prove i): the other cases are similar.

Let $\phi$ be the homomorphism that determines the
$\fgl(\lambda)$-action on $M^{\lambda, s}$.  Then
$\phi(H)=\mathop{\sum}\limits_{i\in\Zee}(s+\eps-2i)E_{ii}$, so
$\phi(e^{tH})=\mathop{\sum}\limits_{i\in\Zee}e^{(s+\eps-2i)t}E_{ii}$
and
$$
\renewcommand{\arraystretch}{1.4}
\begin{array}{l}
\phi((1-e^{-2t})e^{tH})=\phi(e^{tH})-\phi(e^{t(H-2)})
=\mathop{\sum}\limits_{i\in\Zee}e^{(s+\eps-2i)t}(E_{ii}-E_{i-1,
i-1})=\\
\mathop{\sum}\limits_{i\in\Zee}e^{(s-2i)t}
\mathop{\sum}\limits_{j=0}^m\eps^j(E_{ii}-E_{i-1,
i-1})\displaystyle\frac{t^j}{j!}.
\end{array}
$$
Therefore,
$$
\renewcommand{\arraystretch}{1.4}
\begin{array}{l}
(1-e^{-2t})F_\theta(t)=\theta(\phi(e^{tH})-\phi(e^{t(H-2)}))=\\
\mathop{\sum}\limits_{i\in\Zee}e^{(s-2i)t}
\mathop{\sum}\limits_{j=0}^m\displaystyle\frac{\theta_{ij}-\theta_{i-1,
j}}{j!}t^j.\end{array}
$$
Since only finitely many coordinates $\theta_{ij}$ are nonzero,
all the sums are, actually, finite. \end{proof}

\section*{\S 5. Quasi-finite modules over $\fgl(\lambda)$}

In what follows we will show that not all quasi-finite modules
over $\fgl(\lambda)$ can be represented in a canonical form
similar to that of $W_{1+\infty}$ from \cite{KR}.

Let $\fgl^{hol}(\lambda)$ be the holomorphic completion of
$\fgl(\lambda)$, i.e., $\fgl^{hol}_0(\lambda)$ is the algebra of
functions in $H$ holomorphic on the whole complex line $\Cee$ and
$\fgl^{hol}_i(\lambda)=\{f(H)X^i\}$ for $i>0$ whereas
$\fgl^{hol}_i(\lambda)=\{f(H)Y^{-i}\}$ for $i<0$ and $f\in
\fgl^{hol}_0(\lambda)$.  The relations in the completed algebra follow
from (1.1) and (1.4), namely, they are
$$
Xf(H)=f(H-2)X,\; Yf(H)=f(H+2)Y \text{ and }
XY=\frac14(\lambda^2-(H-1)^2).
$$

The following Proposition is proved via the same lines as its
counterpart in \cite{KR}.

\ssbegin{5.1}{Proposition} Let $V$ be quasi-finite
$\fgl(\lambda)$-module.  Then the $\fgl(\lambda)$-action can be
naturally extended to a $\fgl^{hol}_i(\lambda)$-action for $i\neq
0$.
\end{Proposition}

\ssec{5.2} On $\Cee$, introduce an equivalence relation by setting
$$
[s]=\begin{cases} s+2\Zee&\text{ if $s+\lambda\not\in2\Zee+1$}\cr
s-2\Zee_+&\text{ if $s=\pm\lambda+1$}\cr s+2\Zee_+&\text{ if
$s=\pm\lambda-1$}\end{cases}
$$
To each class $[s]$ assign an irreducible $\fsl(2)$-module $M^s$
with a diagonal action of $H$ whose set of eigenvalues coincides
with $[s]$; denote by $\fgl_s(\infty)$ the Lie algebra of linear
transformations of $M^s$ with finitely many nonzero diagonals in
the $H$-diagonal basis.

The following Theorem is similar and is proved via the same
lines as its counterpart in \cite{KR}.

\begin{Theorem} Consider $\mathop{\oplus}\limits_{i=1}^k
M^{s_{i}}_{R_{m_{i}}}$, where the $[s_i]$ are distinct and the
$m_i$ are non-negative integers.  Let $\phi:\fgl(\lambda)\tto\fg=
\mathop{\oplus}\limits_{i=1}^k\fgl_{s_{i}}(\infty; R_{m_{i}})$
determined by this module.  Then, for any quasi-finite
$\fg$-module $V$, any its $\fgl(\lambda)$-submodule is a
$\fg$-submodule.  In particular, if $V$ is irreducible as
$\fg$-module, it is irreducible as a $\fgl(\lambda)$-module.
\end{Theorem}

Following \cite{KR}, let us describe now the structure of
quasi-finite $\fgl(\lambda)$-modules.  Let $\theta$ be the highest
weight of a quasi-finite $\fgl(\lambda)$-module and
$F_\theta(t)=\displaystyle\frac{R(t)}{1-e^{-2t}}$, where
$R(t)=\sum r_i(t)e^{s_{i}t}$ the corresponding formal power
series, see 3.2.  The polynomials $r_i(t)$ will be referred to as
{\it multiplicities}.

\ssec{5.3} For every $s$, denote by $R_s(t)$ the sum of all
quasi-polynomials with the exponents from class $[s]$.  Then
$R(t)=\sum R_s(t)$, where the sum runs over representatives of
different equivalency classes.  Let $R_s(t)=\sum
r_{is}(t)e^{(s-2i)t}$.  The following properties of
quasi-polynomials follow easily from definitions:

i) $R(0)=0$;

ii) the sum of all multiplicities of $R_{\lambda+1}(t)$ as well as
that of $R_{-\lambda+1}(t)$ are equal to 0;

iii) the sum of all multiplicities of $R_{\lambda-1}(t)$ as well as
that of $R_{-\lambda-1}(t)$ are equal to a constant.

Denote: $\Lambda=
[\lambda-1]\cap[\lambda+1]\cap[-\lambda-1]\cap[-\lambda+1]$.  To a
quasi-polynomial $R$ satisfying i)--iii) assign a
$\fgl(\lambda)$-module $V(s)$ as follows:

1) if $s\not\in \Lambda$, let $V(s)$ be the irreducible
$\widehat{\fgl}(\infty)$-module with central charges
$c_j=-\mathop{\sum}\limits_{i\in I}r_i^{(j)}(0)$, where $j=0$, \dots ,
$\max(\deg r_i)$, and the other coordinates of the highest weight are
$\theta_{ij}=\mathop{\sum}\limits_{l\leq
i}(r_l^{(j)}(0)+\delta_{l0}c_j)$;

2) if $s\not\in [\lambda-1]\cap[-\lambda-1]$ and
$R_s=\mathop{\sum}\limits_i r_i(y)e^{s-2i}$, set
$c=-\mathop{\sum}\limits_{i\in I}r_i(0)$ and
$\theta_{ij}=\mathop{\sum}\limits_{l\leq
i}(r_l^{(j)}(0)+\delta_{l0}\delta_{j0}c_j)$; let $V(s)$ be the
corresponding $\fgl^-(\infty)$-module.

3) if $s\not \in [-\lambda+1]\cap[\lambda+1]$, set
$\theta_{ij}=\mathop{\sum}\limits_{l\leq i}r_l^{(j)}(0)$; let $V(s)$
be the corresponding $\fgl^+(\infty)$-module.

\begin{Theorem} Let the quasi-polynomial
$R=\mathop{\sum}\limits_{i=1}^kR_{s_{i}}(t)$ be decomposed with
respect to the distinct equivalence classes of the exponents; let
$R$ satisfy the conditions {\em i)--iii)} above.  Then the
irreducible quasi-finite module with highest weight $R(t)$ is
isomorphic to
$$
V=V(s_1)\otimes \dots\dots V(s_k)\otimes(\alpha\alpha~\tr),
$$
where $\alpha~\tr$ is the irreducible $1$-dimensional module
corresponding to the trace on $\fgl(\lambda)$ and the $V(s_i)$
are constructed according to {\em 1)--3)} above. \end{Theorem}

\begin{proof} By Theorem 4.5 $V$ is irreducible as
$\fgl(\lambda)$-module.  Therefore, it suffices to show that the
restriction of the highest weight of this module onto $\fgl(\lambda)$
is equal to $\displaystyle\frac{R(t)}{1-e^{-2t}}$.  Let $s=s_i$ be one
of the exponents such that $s\not\in\Lambda$.  Then Proposition 5.1
implies that there exists a highest weight $\theta_s$ for
$\widehat{\fgl}(\infty; R_m)$ whose generating function is equal to
$\displaystyle\frac{R_s(t)+c_se^{-(\lambda+1)t}}{1-e^{-2t}}$, where
$c_s=-R_s(0)$.

If $s\in\Lambda$, then the same Proposition implies that the
generating function is equal to
$\displaystyle\frac{R_s(t)}{1-e^{-2t}}$, where the sum of all the
exponents of $R_s$ vanishes.

Finally, the generating function of the trace is equal to
$\alpha\displaystyle\frac{e^{(\lambda-1)t}-e^{-(\lambda+1)t}}{1-e^{-2t}}$.
 Therefore, for our $V$, the generating function is of the form
$\displaystyle\frac{R(t)}{1-e^{-2t}}$, where
$$
\renewcommand{\arraystretch}{1.4}
\begin{array}{l}
R(t)=\mathop{\sum}\limits_{s\not\in\Lambda}(R_s(t)+
c_se^{-(\lambda+1)t})+\\R_{\lambda+1}(t)+
R_{-\lambda+1}(t)+R'{}_{-\lambda-1}(t)+R'{}_{\lambda-1}(t)+
\alpha(e^{(\lambda-1)t}-e^{-(\lambda+1)t})=\\
(\text{by the choice of
}c_s)=\mathop{\sum}\limits_{s\not\in\Lambda}R_s(t)+R_{\lambda+1}(t)+
R_{-\lambda+1}(t)+\\
(R'{}_{-\lambda-1}(t)-\alpha
e^{-(\lambda+1)t})+ (R'{}_{\lambda-1}(t)+
\alpha e^{(\lambda-1)t})=\mathop{\sum}\limits_{s\not\in\Lambda}R_s(t).
\end{array}
$$
\end{proof}

\section*{\S 6. Unitary modules over $\fgl(\lambda)$}

Recall that an {\it anti-involution} of the Lie algebra $\fg$ over
$\Cee$ is an $\Ree$-linear map $\omega:\fg\tto\fg$ such that
$$
\omega(\alpha x)=\bar\alpha\omega(x),\quad
\omega([x, y])=[\omega(x), \omega(y)] \text{ and
}\omega^2=\id
$$
for any $\alpha\in\Cee$ and $x, y\in\fg$.

In presence of an anti-involution $\omega$ we can endow
$V^*$,  the dual of the $\fg$-module $V$, with another
$\fg$-module structure, namely, set
$$
(xl)(v)=l(\omega(x)v)\text{ for any $l\in V^*$, $v\in V$ and
$x\in\fg$}.\eqno{(6.0.1)}
$$
In particular, if $\fg=\fg_-\oplus\fg_0\oplus\fg_+$ and $\omega$
interchanges $\fg_-$ with $\fg_+$, then a $\fg$-homomorphism $M\tto
M^*$ is well-defined, where $M$ is a Verma module with highest weight
vector $v$ and the $\fg$-module structure on $M^*$ is given by
(6.0.1).

Indeed, let $v^*\in M^*$ be such that $v^*(v)=1$ and $v^*(u)=0$ for
all $u$ of lesser weight than that of $v$.  Hence, if $x\in\fg_+$,
then $(xv^*)(u)=v^*(\omega(x)u)=0$, so $v^*$ is also a highest vector
and if $\theta(\omega(H))=\overline{\theta(H)}$, where $\theta$ is the
weight of $v$, then the weight of $v^*$ is also equal to $\theta$.
Hence, there exists a $\fg$-isomorphism $M\tto M^*$ or an Hermitian
$\fg$-invariant form $\langle\cdot, \cdot\rangle$ on $M$.  If
$\langle\cdot, \cdot\rangle$ is positive definite, $M$ is called {\it
unitary}.

In this section we indicate the conditions on the highest weight
of the quasi-finite $\fgl(\lambda)$-module to be unitary.  First,
let us describe automorphisms of $\fgl(\lambda)$.

As is well-known, the description of automorphisms of the algebra
of functions, $\cF$, especially polynomial ones, is a wild
problem. Hence, so is the problem of description of automorphisms
of the Lie algebra $\fder(\cF)$ of the differentiations of $\cF$.
To diminish the amount of automorphisms to a reasonable number, it
is natural to consider the {\it outer} automorphisms only, i.e.,
the classes of automorphisms modulo the group of {\it inner} ones.
Which automorphisms should be considered as inner ones?  In the
case of $\fder(\cF)$ the automorphisms induced by the
automorphisms of $\cF$ naturally qualify.  Similarly, we say that
the automorphism $\phi$ of $\fgl(\lambda)$ is an {\it inner} one
if it also serves as an automorphism of the associative algebra
$\fA_{\lambda}$.

We extend this definition to $LU_\fg(\lambda)$ and call its
automorphism an inner one if it also serves as an automorphism of the
associative algebra $U_\fg(\lambda)$.

\begin{Problem} How to define inner automorphisms for Lie subalgebras of
$LU_\fg(\lambda)$, such as $\fo/\fsp(\lambda)$, to start with? \end{Problem}

Recall that a linear map $\phi: \fg\tto\fg$ is an {\it
anti-automorphism} of the Lie algebra $\fg$ if $\phi([x,
y])=[\phi(y), \phi(x)]$.  For the associative algebras, definition
of an anti-automorphism is similar.  For example, the map $t$ such
that ${}^t|_{\fg}=-\id$ is called the {\it principal
anti-automorphism} of the Lie algebra $\fg$.  The principal
anti-automorphism of $\fg$ can be extended to an anti-automorphism
of the associative algebra $U(\fg)$: ${}^t(x\otimes
y)={}^ty\otimes {}^tx$.  Clearly, $t$ preserves Casimir elements.

\ssbegin{6.1}{Theorem} The group of outer automorphisms of
$\fgl(\lambda)$ is isomorphic to $\Zee/2$ and is generated by
the class of $-t$, where $-t: x\mapsto -{}^tx$.
\end{Theorem}

\begin{proof} Let $\psi\in Aut(\fgl(\lambda))$.  Then $\psi(Y)$,
$\psi(H)$, and $\psi(X)$ span a Lie algebra isomorphic to $\fsl(2)$.
Moreover, they generate $\fA_{\lambda}$ because the weight of
$\psi(X)^n$ is equal to $2n$ with respect to $\psi(H)$.  Hence, there
exists a homomorphism
$$
\tau: U(\fsl(2))\tto \fA_{\lambda}, \quad X\mapsto \psi(X), \;
H\mapsto \psi(H), \; Y\mapsto \psi(Y).
$$
Since $\tau(\Omega)\in\Cee$, there exists a $\mu\in\Cee$ such that
$\tau(\Omega)=\frac12(\mu^2-1)$.  Therefore, $\tau$ is surjective,
hence, an isomorphism $\fA_{\mu}\tto\fA_{\lambda}$.  In a very
difficult technical paper \cite{Di1} Dixmier proved that this can only
happen if $\lambda^2=\mu^2$.  Making use of this, we may assume that
$\lambda=\mu$ and $\tau$ is an automorphism.  Let
$\psi_{1}=\tau^{-1}\psi$.  Then $\psi_{1}|_{\fsl(2)}=\id$.

Therefore, $\psi_{1}$ is an automorphism of $\fgl(\lambda)$, as
$\fsl(2)$-module.  Therefore, $\psi_{1}|_{L^{2i}}=c_{i}\in\Cee$.
Since the Lie algebra $\fgl(\lambda)$ is generated by $L^{2}$ and
$L^{4}$, we deduce that $c_{i}=c_{2}^{i-1}$ for every $i\geq 1$.
Moreover, since $[L^{4}, L^{4}]\supset L^{2}$, it follows that
$c_{2}^{2}= 1$.  Hence, $c_{2}= \pm 1$; if $c_{2}=1$, then
$\psi_{1}=\id$, and if $c_{2}=-1$, then $\psi_{1}=\phi$.  \end{proof}

\ssbegin{6.2}{Theorem} If $\lambda^2\not\in\Ree$, then
$\fgl(\lambda)$ has no real forms, i.e., no involutive anti-linear
automorphisms.
\end{Theorem}

\begin{proof} Let $\omega$ be an involutive anti-linear automorphism of
$\fgl(\lambda)$.  Then $\omega(X)$, $\omega(H)$ and $\omega(Y)$
generate $\fA_{\lambda}$.  Hence, there exists a surjective
homomorphism
$$
\tau: U(\fsl(2))\tto \fA_{\lambda}, \quad \tau=\id \text{ on }X, H, Y;
$$
hence, we obtain an isomorphism $\tau: \fA_{\mu}\tto \fA_{\lambda}$.
By \cite{Di1} this may only happen if $\lambda^2=\mu^2$.  Hence, we
may assume that $\lambda=\mu$ and $\tau$ is an automorphism.  Then
$\omega_{1}=\tau^{-1}\omega$ is an antilinear automorphism of
$\fgl(\lambda)$ such that
$$
\omega_{1}(X)=X, \; \omega_{1}(H)=H, \; \omega_{1}(Y)=Y.
$$
Set $z_{i}=(\ad~X)^i (Y^2)$; in particular, $z_0=z$.  Then
$\Span(z_{i}: i=1, \dots , 4)$ form a basis of $L^4$ and if
$\omega_{1}(Y^2)=cY^2$, then, clearly, $\omega_{1}(z_{i})=cz_{i}$
for $i=1$, \dots, 4.  According to \cite{GL}, in $\fgl(\lambda)$
the following relations hold:
$$
\renewcommand{\arraystretch}{1.4}
\begin{array}{l}
    3[z_{1}, z_2]-2[z, z_3]=24(\lambda^{2}-4)Y\\
    4[z_3, [z, z_1]]-3[z_2, [z, z_2]]=576(\lambda^2-9)z
\end{array}
$$
Having applied $\omega_{1}$ to both sides of these relations we get
$$
\renewcommand{\arraystretch}{1.4}
\begin{array}{l}
    c^2(\lambda^{2}-4)=\bar\lambda^{2}-4\\
    c^3(\lambda^2-9)=c(\bar\lambda^2-9)
\end{array}
$$
Therefore, $c^2=1$ and $\lambda^2=\bar\lambda^2$. \end{proof}

\begin{Corollary} If $\lambda^2\neq\bar\lambda^2$, then
$\fgl(\lambda)$ has no involutive anti-linear automorphisms.
\end{Corollary}

\begin{proof} Let $\omega$ be an involutive anti-linear automorphism;
clearly $\tau: X\longleftrightarrow Y$ and $\tau(H)=H$ determines
an anti-automorphism of $\fgl(\lambda)$.  Then $\tau\circ\omega$
is an involutive anti-linear automorphism.
\end{proof}

\ssec{6.3} If $\lambda^2=\bar\lambda^2$, then $\fA_{\lambda}$
possesses an involutive anti-linear automorphism $\omega$ given on
generators by the formulas:
$$
\omega(X)=Y, \quad \omega(Y)=X, \quad
\omega(H)=H\eqno{(6.3)}
$$
In what follows, the unitary modules are considered with
respect to this automorphism.

\begin{Theorem} Let $\lambda^ 2$ be real and let
$F_\theta(t)=\displaystyle\frac{R(t)}{1-e^{-2t}}$, where $R(t)$ is
a quasi-polynomial none of whose exponents $s_i$ belong to
$\Lambda=[\lambda-1]\cap[\lambda+1]\cap[-\lambda-1]\cap
[-\lambda+1]$. The $\fgl(\lambda)$-module with character
$F_\theta(t)$ is unitary if and only if
$$
F_\theta(t)=\sum
n_i\frac{e^{-(\lambda+1)t}-e^{s_it}}{1-e^{-2t}},\text{ where
}n_i\in\Zee_+ .
$$ \end{Theorem}

\begin{proof} Let $V$ be the irreducible module with character
$F_\theta(t)$ and $P(H)$ the annihilator of $Yv$, where $v$ is the
highest weight vector.  Then $P(H)$ is the characteristic
polynomial of the operator $H$ in $V_{-1}$.  But since the unitary
form is invariant and $\omega(H)=H$, it follows that $H$ is
self-adjoint; hence, all the roots of $P$ are real.  Further, if
$\alpha$ is a multiple root of multiplicity $m>1$, the polynomial
$P$ is of the form $P=(H-\alpha)^mQ(H)$.  For $u=
(H-\alpha)^{m-1}Q(H)v$, we have
$$
\langle u, u\rangle=\langle (H-\alpha)^{m-1}Q(H)v,
(H-\alpha)^{m-1}Q(H)v\rangle=\langle Q(H)v,
(H-\alpha)^{2m-2}Q(H)v\rangle=0.
$$
Thanks to Hermitian property, $u=0$. Thus, $m=1$.

By Theorem 5.3 $V$ is of the form $V(s_1)\otimes\dots\otimes
V(s_k)$, where the $V(s_i)$ are irreducible
$\widehat{\fgl}(\infty)$-modules. Clearly, $V$ is unitary if and
only if each $V(s_i)$ is unitary.  Now, thanks to \cite{KR} we
know that the $\widehat{\fgl}(\infty)$-module $V(s)$ is unitary if
and only if $\theta_i-\theta_{i+1}+\delta_{i0}c\in\Zee_+$ for all
$i$, where the $\theta_i$ are the coordinates of the highest
weight and $c$ the value of the central charge.  \end{proof}

\section*{\S 7.  $\fgl(\lambda)$ and the symmetric group}
In this section we deduce an explicit realization of certain
irreducible $\fgl(\lambda)$-modules.  Namely, we decompose the tensor
powers of the Verma module over $\fsl(2)$, indicate the corresponding
characteristic polynomials and $q$-characters.

\ssbegin{7.1}{Lemma} Let $A$ be an associative
algebra; let $\fS_n$ naturally act on $A^{\otimes
n}$. Then the algebra $(A^{\otimes
n})^{\fS_n}$ of $\fS_n$-invariants is generated by the
elements of the form
$$
a\otimes 1\otimes 1\otimes \dots\otimes 1+\dots
1\otimes 1\otimes 1\otimes \dots\otimes a.
$$
\end{Lemma}

\begin{proof} Denote: $s(a_1, \dots,
a_n)=\mathop{\sum}\limits_{\sigma\in\fS_{n}}a_{\sigma(1)}\otimes
\dots\otimes a_{\sigma(n)}$.  Let $B$ be the algebra generated by the
$s(a, 1, \dots, 1)$ for all $a\in A$.  Let $|s(a_1, \dots, a_n)|$ be
the number of the $a_i$ distinct from 1.  Let us prove that $B\simeq
(A^{\otimes n})^{\fS_n}$.

Let us carry on an induction on $|s(a_1,\dots, a_n)|$ in order to
prove that $s(a_1, \dots, a_n)\in B$.  Indeed, if $|s(a_1, \dots,
a_n)|=1$, by definition $s(a_1, \dots, a_n)=s(a, 1, \dots, 1)\in B$.

Let $|s(a_1,\dots, a_n)|=l>1$. Consider
$$
s(a_1,\dots, a_{l-1}, 1, \dots, 1)s(a_l, 1, \dots, 1)=\alpha
s(a_1,\dots, a_{l}, 1, \dots, 1)+\dots,
$$
where $\alpha$ is a nonzero constant and the dots stand for the linear
combination of the terms $s(b_1,\dots, b_{m}, 1, \dots, 1)$
with $m<l$. \end{proof}

\ssbegin{7.2}{Theorem} Let $V=M^{\lambda-1}$ be
the Verma module with highest weight
$\lambda-1$ over $\fsl(2)$ (hence, a $\fgl(\lambda)$-module). Then
$$
V^{\otimes n}=\mathop{\oplus}\limits_\nu V^{\nu}\otimes S^{\nu},
$$
where $V^{\nu}$ is an irreducible $\fgl(\lambda)$-module and
$S^{\nu}$ is an irreducible $\fS_n$-module and $\nu$ runs over the
partitions of $n$.
\end{Theorem}

\begin{proof} Clearly, $V$ is irreducible not only as
$\fgl(\lambda)$-module but also as $\fA_\lambda$-module.  Then,
obviously, $W= V^{\otimes n}$ is irreducible as $\fA_\lambda^{\otimes
n}$-module.

The image of $U(\fgl(\lambda))$ in $\End(W)$ coincides with the
subalgebra generated by $s(a, 1, \dots, 1)$ for
$a\in\fgl(\lambda)$, so by Lemma 7.1 it is isomorphic to
$(\fA_\lambda^{\otimes n})^{\fS_{n}}$.

Let us decompose $W$ into isotypical $\fS_{n}$-modules:
$W=\mathop{\oplus}\limits_\nu W^{\nu}$ and represent each $
W^{\nu}$ in the form $ W^{\nu}=V^{\nu}\otimes S^{\nu}$, where
$V^{\nu}$ is an irreducible $\fgl(\lambda)$-module.  There is no
unique way to do this, for example, set $V^{\nu}=e_\nu(W)$ for any
minimal idempotent in $\Cee[\fS_n]$ corresponding to the partition
$\nu$.

To show that $V^{\nu}$ is irreducible as a $\fgl(\lambda)$-module,
consider $V^{\nu}_1=\Hom_{\fS_{n}}(S^\nu, W)$.  As a
$\fgl(\lambda)$-module, $V^{\nu}_1$ is isomorphic to $V^{\nu}$.
Let $\phi, \psi\in V^{\nu}_1$ and $\phi\neq 0$.  Let us show that
there exists a $u\in U(\fgl(\lambda))$ such that $u\phi=\psi$.
Indeed, since $W$ is irreducible as $\fA_\lambda^{\otimes
n}$-module, then the density theorem (\cite{L}, Ch.  XVII, \S 3,
Th.1) states that there is a $w\in \fA_\lambda^{\otimes n}$ such
that $w\phi(v_i)=\psi(v_i)$ , where the $v_i$ form a basis of
$S^\nu$.  (Since $S^\nu$ is irreducible and $\phi\neq 0$, then the
vectors $\psi(v_i)$ are linearly independent for $i=1$, \dots ,
$\dim S^\nu$.)

Let us average the element $w\phi$ with respect to $\fS_{n}$:
$$
(w\phi)^{\sharp}=
\frac{1}{|\fS_{n}|}\sum \sigma(w\phi)\sigma^{-1}=
\frac{1}{|\fS_{n}|}\sum
\sigma(w)\sigma^{-1}(\phi)=w^{\sharp}\phi.
$$
But, on the other hand, since $w\phi=\psi$, we deduce that
$(w\phi)^{\sharp}=(\psi)^{\sharp}=\psi$, i.e., $w^{\sharp}\phi=\psi$
and we may assume that $w\in (\fA_\lambda^{\otimes n})^{\fS_{n}}$.
So, there exists a $u\in U(\fgl(\lambda))$ such that $u\phi=\psi$.  In
other words, $V_1^\nu$ is irreducible.
\end{proof}

\ssbegin{7.2.1}{Corollary} {\em 1)} Let $\nu=(\nu_1, \dots , \nu_n)$.
The generating function corresponding to the highest weight $\nu$ is
of the form
$$
\mathop{\sum}\limits_{i=1}^n \nu_ie^{(\lambda -2i+1)t}.\eqno{(7.2.1.1)}
$$

{\em 2)} Represent $\nu$ in the form $\nu=(\theta_1^{\alpha_{1}}\dots
\theta_m^{\alpha_{m}})$, where $\theta_1>\dots>\theta_m>0$ and
$\alpha_{i}\neq 0$ for all $i$ and where $\theta^{\alpha}$ denotes the
product of $\alpha$ copies: $\theta\dots \theta$.  Then the
characteristic polynomial of $V^\nu$ is equal to
$$
P(H)=\prod_{i=1}^m(H-\lambda-2\alpha_1-2\alpha_2-\dots
-2\alpha_i-1).
$$
\end{Corollary}

\begin{proof} Heading 1) follows from Theorem 4.5. Multiply
(7.2.1.1) by $1-e^{-2t}$; we see that the product is of the form
$$
\nu_1e^{(\lambda -1)t}+\mathop{\sum}\limits_{i=1}^{n-1}
(\nu_{i+1}-\nu_i)e^{(\lambda -2i+1)t}-\nu_ne^{(\lambda-2n
-1)t}.\eqno{(7.2.1.2)}
$$
Therefore, sum (7.2.1.2) without the first summand is a solution of an
ordinary differential equation with constant coefficients whose
characteristic equation is precisely of the form indicated.
\end{proof}

\ssbegin{7.2.2}{Corollary} Set $a=e^\lambda$, $q=e^{-\alpha}$, where
$\alpha$ is the positive root of $\fsl(2)$.  Then the $q$-character of
$V^\nu$ is equal to
$$
\chi_\nu=\frac{a^{|\nu|}q^{n(\nu)}}
{\prod_{x\in\nu}(1-q^{h(x)})},\eqno{(7.2.2)}
$$
where $|\nu|$ is the number of cells in the Young tableau
corresponding to $\nu$, $n(\nu)=\mathop{\sum}\limits_{i\geq 1}\nu_i$
and $h(x)$ is the length of the hook corresponding to the cell $x$.
\end{Corollary}

\begin{proof} By \cite{Mac}, Example 2 in \S 3, for the $S$-function
$s_\nu$ we have
$$
\chi_\nu=s_\nu(a, aq, aq^2, \dots, )=a^{|\nu|}s_\nu(1, q,q^2,
\dots, ) =\text{rhs of }(7.2.2).
$$
\end{proof}

\section*{\S 8. Orthogonal polynomials for $\fgl(n)$}
For an overview, see \cite{NU}, \cite{NSU}.  (Setting
$T_{1}(\alpha_{i})\dots T_{l}(\alpha_{i})=1$ for $l=0$ we make
formulas (8.1.1) and the like look uniform.)

\ssbegin{8.1}{Theorem} In $\fgl(n)$:

{\em i)} Consider the basis $e_{kl}=(\ad ~Y)^{k-l}(X^k)$ for $0\leq
k\leq n-1$ and $-k\leq l\leq k$.  We have
$$
\langle e_{kl}, e_{k'l'}\rangle=\delta_{k, k'}\delta_{l+l', 0}.
$$

{\em ii)} Determine the elements $f_{kl}$ from the equations ($0\leq
l\leq k$)
$$
(\ad ~Y)^{k-l}(X^k)=X^lf_{kl}\text{ and }
(\ad ~Y)^{k+l}(X^k)=f_{k, -l}Y^l.
$$
Set $T_{i}(H)=\frac14(n^{2}-(H+2i-1)^2)$ and $\alpha_{i}=n-2i+1$
for $i=1$, \dots , $n$.  For a fixed $l\geq 0$ and any $k\geq l$,
the polynomials $f_{kl}$ form an orthogonal basis with respect to
the form
$$
\langle f, g\rangle=\begin{cases}
\mathop{\sum}\limits_{i=l+1}^nf(\alpha_{i})
g(\alpha_{i})T_{1}(\alpha_{i})\dots T_{l}(\alpha_{i})&\text{ for
$l>0$}\cr \mathop{\sum}\limits_{i=1}^nf(\alpha_{i})
g(\alpha_{i})&\text{ for $l=0$}.\end{cases}\eqno{(8.1.1)}
$$

{\em iii)} Up to a constant factor the polynomials $f_{kl}$ coincide
with the Hahn polynomial of one discrete variable:
$$
f_{kl}(H)=\, {}_{3}F_{2}\left(\begin{matrix}l-k, l+k+1,
\frac12(1-n-H)\cr l+1,
l+1-n\end{matrix}~\vert~ 1\right)\times T_{0}(\alpha_{l+1})\dots
T_{0}(\alpha_{k}),\eqno{(8.1.2)}
$$
where
$$
{}_{3}F_{2}\left(\begin{matrix}\alpha_{1}, \alpha_{2}, \alpha_{3}\cr
\beta_{1}, \beta_{2}\end{matrix}\mid
z\right)=\mathop{\sum}\limits_{i=0}^\infty
\displaystyle\frac{(\alpha_{1})_{i}(\alpha_{2})_{i}
(\alpha_{3})_{i}}{(\beta_{1})_{i}(\beta_{2})_{i}}\,
\displaystyle\frac{z^i}{i!}\eqno{(8.1.3)}
$$
is a generalized hypergeometric function, $(\alpha)_{0}=1$  and
$(\alpha)_{i}=\alpha(\alpha+1)\dots(\alpha+i-1)$ for $i>0$.
\end{Theorem}

\begin{proof} i) First, observe that the subspaces $L^{2k}$ in the
decomposition of $\fgl(n)=L^{0}\oplus L^{2}\oplus\dots \oplus
L^{2n-2}$ are pairwise orthogonal.  Indeed, the form $A, B\mapsto
\tr~AB$ defines an invariant pairing of $L^{k}$ and $L^{l}$, hence,
an $\fsl(2)$-homomorphism $L^{k}\tto L^{l}$ which is only possible
if $k=l$. Moreover, it is clear that $\langle \fgl(n)_{k},
\fgl(n)_{l}\rangle\neq 0$ if and only if $k+l=0$. Now observe that
$\fgl(n)_{l}\cap L^{2k}=\Span((\ad Y)^{k-l}(X))$. This proves i).

ii) Let $X^lf(H)\in \fgl(n)_{l}$ and $g(H)Y^l\in \fgl(n)_{-l}$. Then
$$
\tr~(g(H)Y^lX^lf(H))=\tr~(f(H)g(H)T_{1}(H)\dots T_{l}(H)).
$$
As is easy to verify
$$
\tr~f(H)=\mathop{\sum}\limits_{1\leq i\leq n}f(\alpha_{i})\text{ for any
}f(H)\in\fgl(n)_{0}.
$$
This implies (8.1.1).

Moreover, i) implies that for any fixed $l$ lying between $-k$ and
$k$ the polynomials $f_{kl}$ and $f_{k, -l}$ form two mutually
dual bases, i.e., $\langle f_{kl}, f_{kl'}\rangle =\delta_{l+l', 0}$.
But the degrees of polynomials $f_{kl}$ and $f_{k, -l}$ are equal,
which means that they coincide up to a constant factor. This proves ii).

iii) follows from the fact that orthogonal polynomials are uniquely
determined by the weight function and the interval over which we
consider the scalar product. Eq. (8.1.2) follows from comparison of the
coefficients of the leading terms. \end{proof}

\begin{Remark} 1) Hahn polynomials are determined for three parameters
$\alpha$, $\beta$ and $N$ as
$$
h_{p}^{(\alpha, \beta)}(z, N)=
\frac{(-1)^p\Gamma(N)(\beta)_{p}}{p!\Gamma(N-p)} \;
{}_{3}F_{2}\left(\begin{matrix}-p, \alpha+\beta+p+1, -z\cr
\beta+1, l-N\end{matrix}~\vert~
1\right ).  \eqno{(8.1.4)}
$$
Our $f_{kl}$ coincides, up to a factor, with $h_{p}^{(\alpha,
\beta)}(z, N)$ at $\alpha=\beta=l$, $p=k-l$, $z= \frac12(H+n-1)$
and $N=n-1$.

2) Having finished this paper we have realized that sometimes it is
more convenient to consider the form $A, B\mapsto \tr~AB^t$, i.e., an
invariant pairing on each $\fgl(n)_{k}$, cf. \cite{S}.
\end{Remark}

\ssec{8.2} Let us show now how the main properties of the Hahn
polynomials follow from the representation theory of $\fsl(2)$.
For any $f(H)\in \Cee[H]$, set
$$
\triangle f(H)=f(H+2)-f(H), \quad \nabla f(H)=f(H)-f(H-2).\eqno{(8.1.5)}
$$

\begin{Theorem} Consider $\fgl(n)$ as $\fsl(2)$-module with
respect to the image of the principal embedding, let
$\Omega$
be the quadratic Casimir operator $(1.2)$.

{\em i)} $\Omega$ is self-adjoint with respect to the form
$\langle\cdot, \cdot\rangle$ and the polynomials $X^lf_{kl}$ are
eigenfunctions of $\Omega$ corresponding to eigenvalue $2k(k+1)$.  The
polynomials $f_{kl}$ satisfy the difference equation
$$
T_0(H)\nabla \triangle (f)-(l+1)(H+l)\triangle
(f)+(k-l)(k+l+1)f=0.\eqno{(8.2.i))}
$$

{\em ii)}
$f_{kl}=\begin{cases}\displaystyle\frac{\nabla^{k-l}(T_{1}\dots
T_{k})}{T_{1}\dots T_{l}}&\text{if $l>0$}\cr \nabla^{k}(T_{1}\dots
T_{k})&\text{if $l=0$.}\end{cases}$

{\em iii)} $\langle f_{kl}, f_{kl}\rangle=\displaystyle
\frac{(k-l)!(k!)^2}{(k+l)!(2k+1)}
n(n^2-1^{2})\dots (n^2-k^{2})$.\end{Theorem}

\begin{proof} i) Since the form $\langle\cdot, \cdot\rangle$ is
$\fsl(2)$-invariant, $\langle[w, u], v\rangle=-\langle u, [w,
v]\rangle$ for any $u, v\in\fgl(n)$ and $w\in\fsl(2)$.  Denote the
$U(\fsl(2))$-action in $\fgl(n)$ by $*$; by induction we easily
deduce that $\langle w*u, v\rangle=-\langle u, w^t*v\rangle$ where
now $w\in U(\fsl(2))$ and $t$ is the principal anti-involution.
Since, as is easy to verify, $\Omega^t=\Omega$, we have $\langle
\Omega*u, v\rangle=\langle u, \Omega*v\rangle$, i.e., $\Omega$ is
self-adjoint. Since $\fgl(n)=\oplus L^{2i}$ and $X^lf_{kl}=(\ad
Y)^{k-l}(X^k)\in L^{2k}$, it follows that $X^lf_{kl}$ is an
eigenfunction of $\Omega$ with eigenvalue $2k(k+1)$.  Applying
$\Omega$ to $X^lf_{kl}$ we obtain eq.  (8.2.i).

ii) Recall the identity
$$
(\ad y)^p(a)=\mathop{\sum}\limits_{j=0}^p(-1)^j\binom{p}{j}y^{p-j}ay^j.
$$
Set $y=Y$, $a=X^k$, $p=k-l$ and multiply by $Y^l$ from the left. We get
$$
\renewcommand{\arraystretch}{1.4}
\begin{array}{l}
Y^lX^lf_{kl}=Y^l(\ad Y)^{k-l}(X^k)=
Y^l\left(\mathop{\sum}\limits_{j=0}^{k-l}(-1)^j\binom{k-l}{j}
Y^{k-l-j}X^kY^j\right)=\\
\mathop{\sum}\limits_{j=0}^{k-l}(-1)^j\binom{k-l}{j}
y^{k-j}X^kY^j=\mathop{\sum}\limits_{j=0}^{k}(-1)^j\binom{k}{j}
T_{j-1}\dots T_{0}T_{1}\dots T_{k-j}. \end{array}
$$
But $\nabla^{k-l}(f)(H)=\mathop{\sum}\limits_{j=0}^{k}(-1)^j\binom{k}{j}
f(H-2j)$ for any function $f$, so
$$
Y^lX^lf_{kl}=\nabla^{k-l}(T_{1}\dots T_{k}).
$$
Since $Y^lX^l=T_{1}\dots T_{l}$, we are done.

iii) The module $L^{2k}$ is the linear span of the vectors $v_{l}=(\ad
Y)^{k-l}(X^k)$ for $-k\leq l\leq k$; hence,
$$
\langle v_{l}, v_{m}\rangle=-\langle v_{l-1},
v_{m+1}\rangle=\dots=(-1)^l\langle v_{0}, v_{l+m}\rangle
$$
by invariance of the form. Hence,
$\langle v_{l}, v_{-l}\rangle=(-1)^l\langle v_{0}, v_{0}\rangle$.
But $v_{l}=X^lf_{kl}=(\ad ~Y)^{k-l}(X^k)$ and $v_{-l}=f_{k, -l}Y^l=
(\ad ~Y)^{k+l}(X^k)$.

Now recall that $f_{kl}$ and $f_{k, -l}$ are identical up to a
constant factor and compare the leading coefficients. We see that
$$
f_{k, -l}=(-1)^l\frac{(k+l)!}{(k-l)!}f_{k, l}.
$$
Hence,
$$
\langle f_{k, l}, f_{k, l}\rangle=\tr~(f_{k, l}Y^l\cdot X^lf_{k,
l})=(-1)^l\displaystyle\frac{(k+l)!}{(k-l)!}\langle v_{l},
v_{-l}\rangle= \displaystyle\frac{(k+l)!}{(k-l)!}\langle v_{0},
v_{0}\rangle.
$$
It remains to demonstrate that
$$
\langle v_{0}, v_{0}\rangle=\displaystyle\frac{n(n^2-1^{2})\dots
(n^2-k^{2})}{2k+1}(k!)^2.
$$

Obviously, $f_{k}=\left((\ad ~Y)^{k}(X^k)\right)^{2}\in U(\fsl(2))$.
So $\langle v_{0}, v_{0}\rangle=\tr_{n}\phi_{n}(f_{k})$, where
$\tr_{n}$ is the trace on $\fgl(n)$ and $\phi_{n}$ is the homomorphism
of $U(\fsl(2))$ into $\fgl(n)$ induced by the principal embedding.

It is clear now that $P_k(n)=\tr_{n}(\phi_{n}(f_{k}))$ is a polynomial of
degree $2k+1$; moreover, this polynomial is an odd one. But
$\phi_{n}(X^k)=0$ for $n\leq k$ and $\phi_{n}(f_{k})=0$; hence,
$P_k(n)=0$ if $n\leq k$ and $P_k(-n)=0$ because $P_k$ is odd.
Therefore,
$$
P_k(n)=c_{k}n(n^2-1^{2})\dots
(n^2-k^{2}).
$$
To calculate the constant $c_{k}$, it suffices to compute
$P_k(k+1)=\tr_{k+1}(\phi_{k+1}(f_{k}))$. But in this case
$X^k=(k!)^{2}e_{1, k+1}$ and
$$
(\ad ~Y)^{k}(X^k)=(k!)^{2}(\ad ~Y)^{k}(e_{1, k+1})=(k!)^{2}
\mathop{\sum}\limits_{j=0}^{k}(-1)^j\binom{k}{j}e_{k+1-j, k+1-j}.
$$
Thus,
$$
\tr_{k+1}(\phi_{k+1}(f_{k}))=(k!)^{4}\mathop{\sum}\limits_{j=0}^{k}
\binom{k}{j}^2=
(k!)^{4}\binom{2k}{k}
$$
implying $c_{k}=\displaystyle\frac{(k!)^{2}}{2k+1}$. \end{proof}

\ssbegin{8.3.1}{Proposition} The following
orthogonality relations {\em (i)} and {\em (ii)} hold:
$$
\mathop{\sum}\limits_{i} f_{kl}(\alpha_{i})f_{k_{1},
l}(\alpha_{i})T_{1}(\alpha_{i})\dots T_{l}(\alpha_{i})=\delta_{k,
k_{1}} c_{k, l}, \leqno{(i)}
$$
where $c_{k,
l}=\displaystyle\frac{(k-l)!(k!)^2}{(k+l)!(2k+1)}n(n^2-1^{2})\dots
(n^2-k^{2})$.
$$
\mathop{\sum}\limits_{k=l}^{n-1}\frac{1}{c_{kl}}f_{kl}(\alpha_{i})
f_{kl}(\alpha_{j})
T_{1}(\alpha_{i})\dots T_{l}(\alpha_{i})=\delta_{ij}\text{ for any $i,
j>l$}.\leqno{(ii)}
$$
\end{Proposition}

\begin{proof} (i) follows from Theorem 8.2, iii).  To prove ii),
express $e_{ij}$ via $X^lf_{kl}$ and $f_{k, -l}Y^l$.  We have, in
particular, $e_{i, i+l}=\mathop{\sum}
\limits_{k=l}^{n-1}\alpha_{ik}X^lf_{kl}$ for any $k$ between $l$ and
$n-1$.  But, as is easy to verify,
$$
\renewcommand{\arraystretch}{1.4}
\begin{array}{l}
    f_{kl}Y^l=(\mathop{\sum}
\limits_{i=1}^{n} f_{kl}(\alpha_{i})e_{ii})(\mathop{\sum}
\limits_{i=1}^{n-1}e_{i+1, i})^l=\\
(\mathop{\sum}
\limits_{i=1}^{n} f_{kl}(\alpha_{i})e_{ii})(\mathop{\sum}
\limits_{i=1}^{n-l}e_{i+l, i})=\mathop{\sum}
\limits_{i=1}^{n-l} f_{kl}(\alpha_{i+l})e_{i+l, i}.
\end{array}
$$
Hence, by Theorem 8.2 we have
$$
f_{kl}(\alpha_{i+l})=\langle e_{i, i+l}, f_{kl}Y^l\rangle= \alpha_{ik}
\langle f_{kl}Y^l, X^lf_{kl}\rangle=\alpha_{ik}c_{kl}.
$$
In other words,
$$
e_{i, i+l}=\mathop{\sum}
\limits_{k=l}^{n-1}\frac{f_{kl}(\alpha_{i+l})}{c_{kl}}X^lf_{kl}.
$$
Similarly,
$$
e_{i+l, i}=\mathop{\sum}
\limits_{k=l}^{n-1}\frac{f_{kl}(\alpha_{i+l})T_{1}(\alpha_{i+1})\dots
T_{l}(\alpha_{i+l-1})}{c_{kl}}f_{kl}Y^l.
$$
Therefore, by setting $i'=i+l$, $j'=j+l$  we deduce that
$$
\renewcommand{\arraystretch}{1.4}
\begin{array}{l}
  \langle e_{i+l, i}, e_{j, j+l}\rangle=\mathop{\sum}
\limits_{i=l}^{n-1}
\frac{1}{c_{kl}}f_{kl}(\alpha_{i+l})f_{kl}(\alpha_{j+l})
T_{1}(\alpha_{i+1})\dots
T_{l}(\alpha_{i+l-1})=\\
\mathop{\sum}
\limits_{i=l}^{n-1} \frac{1}{c_{kl}}f_{kl}(\alpha_{i'})f_{kl}(\alpha_{j'})
T_{1}(\alpha_{i'-1})\dots
T_{l}(\alpha_{i'-l+1})=\\
\mathop{\sum}
\limits_{i=l}^{n-1} \frac{1}{c_{kl}}f_{kl}(\alpha_{i'})f_{kl}(\alpha_{j'})
T_{1}(\alpha_{i'})\dots
T_{l}(\alpha_{i'}).
\end{array}
$$
\end{proof}

\ssec{8.3.2} Expressing Casimir elements $\Omega_{2}$,
$\Omega_{3}$, \dots , $\Omega_{n}$ for $\fgl(n)$ in terms of the
orthogonal polynomials rather than matrix units, we can derive
various identities relating the polynomials.  Some of these
identities might be even new, at least, for non-integer values of
$n$, cf.  9.3.  For example, for $\Omega_{2}$, we have:

\begin{Proposition} In $\fgl(n)$, we have
$$
\mathop{\sum}\limits_{k=0}^{n-1}\left(
\frac{1}{c_{k0}}f_{k0}^2+2\mathop{\sum}\limits_{l=1}^{k}\frac{1}{c_{kl}}
f_{kl}^2
T_{1}\dots T_{l}-\frac{2k+1}{n}\right)=-H.
$$
\end{Proposition}

\begin{proof} Let us rewrite the Casimir operator for $\fgl(n)$
expressed in terms of matrix units
$$
\renewcommand{\arraystretch}{1.4}
\begin{array}{l}
    \Omega_{2}=\mathop{\sum}\limits_{k=1}^{n}e_{kk}\otimes e_{kk}+
    \mathop{\sum}\limits_{i<j}(e_{ij}\otimes e_{ji}+e_{ji}\otimes e_{ij})=\\
\mathop{\sum}\limits_{k=1}^{n}e_{kk}\otimes e_{kk}+
\mathop{\sum}\limits_{k=1}^{n}(n-2k+1)e_{kk}+
    2\mathop{\sum}\limits_{i<j}e_{ji}\otimes e_{ij})
\end{array}
$$
via orthogonal polynomials:
$$
\renewcommand{\arraystretch}{1.4}
\begin{array}{l}
    \Omega_{2}=\mathop{\sum}\limits_{k=0}^{n-1}
\frac{1}{c_{k0}}f_{k0}\otimes f_{k0}+\mathop{\sum}\limits_{l=1}^{n-1}
\mathop{\sum}\limits_{k=l}^{n-1}\frac{1}{c_{kl}}(X^lf_{kl}\otimes
f_{kl}Y^l+
f_{kl}Y^l\otimes X^lf_{kl})=\\
\mathop{\sum}\limits_{k=0}^{n-1}
\frac{1}{c_{k0}}f_{k0}\otimes f_{k0}+2\mathop{\sum}\limits_{l=1}^{n-1}
\mathop{\sum}\limits_{k=l}^{n-1}\frac{1}{c_{kl}}
f_{kl}Y^l\otimes X^lf_{kl}+\mathop{\sum}\limits_{k>l}
\frac{1}{c_{kl}}[X^lf_{kl},  f_{kl}Y^l].
\end{array}
$$
Let us now transform the last sum. For this, consider a
homomorphism $\varphi: U(\fgl(n))\tto \fgl(n)$ induced by the identity
representation. Since
$$
\frac{1}{c_{k0}}f_{k0}\otimes f_{k0}+\mathop{\sum}\limits_{l=1}^{k}
\frac{1}{c_{kl}}(X^lf_{kl}\otimes
f_{kl}Y^l+
f_{kl}Y^l\otimes X^lf_{kl})
$$
is
$\fsl(2)$-invariant, its image in $\fgl(n)$ is also an
$\fsl(2)$-invariant, i.e., is a constant, or better say, scalar
matrix, $D(k)$:
$$
\frac{1}{c_{k0}}f_{k0}\cdot f_{k0}+2\mathop{\sum}\limits_{l=1}^{k}
\frac{1}{c_{kl}}
f_{kl}Y^l\cdot X^lf_{kl}+\mathop{\sum}\limits_{l=1}^k
\frac{1}{c_{kl}}[X^lf_{kl},  f_{kl}Y^l]=D(k).
$$

Let us calculate the trace of both sides. We see that
$$
nD(k)=\frac{1}{c_{k0}}\tr(f_{k0}^2)+2\mathop{\sum}\limits_{l=1}^{k}
\frac{1}{c_{kl}}\tr(f_{kl}^2T_{1}\dots T_{l})= 2k+1
$$
implying $D(k)=\displaystyle\frac{2k+1}{n}$. Therefore,
$$
\mathop{\sum}\limits_{l=1}^k \frac{1}{c_{kl}}[X^lf_{kl},
f_{kl}Y^l]=\displaystyle\frac{2k+1}{n}-\left( \frac{1}{c_{k0}}f_{k0}\cdot
f_{k0}+2\mathop{\sum}\limits_{l=1}^{k} \frac{1}{c_{kl}} f_{kl}Y^l\cdot
X^lf_{kl}\right).
$$
Summing over $k$ we obtain
$$
\mathop{\sum}\limits_{k=0}^{n-1}\left( \frac{1}{c_{k0}}f_{k0}\cdot
f_{k0}+2\mathop{\sum}\limits_{l=1}^{k} \frac{1}{c_{kl}} f_{kl}Y^l\cdot
X^lf_{kl}-\displaystyle\frac{2k+1}{n}\right) =-\mathop{\sum}\limits_{k=0}^{n-1}\mathop{\sum}\limits_{l=1}^k
\frac{1}{c_{kl}}[X^lf_{kl}, f_{kl}Y^l].
$$
The linear parts of the same Casimir operator but expressed in
different bases coincide, so the last
sum is equal to $-\mathop{\sum}\limits_{i=1}^{n}(n-2i+1)e_{ii}=-H$.
\end{proof}

\section*{\S 9.  Orthogonal polynomials and $\fgl(\lambda)$}

The above results for $\fgl(n)$, $n=1$, $2$, etc.  also hold {\it mutatis
mutandis} for $\fgl(\lambda)$ with any complex $\lambda\not\in
\Zee\setminus \{0\}$.  In this section we only consider such values of
$\lambda$.

Observe that $\fgl(\lambda)$, unlike $\fgl(\infty)$ or
$\fgl^{\pm}(\infty)$, has no basis consisting of matrix units. But we
can represent it in the form
$$
\fgl(\lambda)=\oplus^\infty_{k=0}L^{2k}\; \text{ and }\;  \fgl(\lambda)=
\oplus^\infty_{i=-\infty}\fgl(\lambda)_{i},
$$
where $\fgl(\lambda)_{0}=\Cee[H]$ while
$\fgl(\lambda)_{i}=X^i\fgl(\lambda)_{0}$ and
$\fgl(\lambda)_{-i}=\fgl(\lambda)_{0}Y^i$ for $i>0$.

As was shown above, there exists a trace on $\fgl(\lambda)$ whose
restriction onto $\fgl(\lambda)_{0}$ has for generating function
$$
\frac{e^{\lambda t}-e^{-\lambda t}}{e^{t}-e^{-t}}\text{ if }\lambda
\neq 0 \text{ and }\frac{2t}{e^{t}-e^{-t}}\text{ if }\lambda= 0.
$$
In the first case we have normalized the trace so that
$\tr(1)=\lambda$ by analogy with the finite dimensional case when
the scalars are naturally represented by scalar matrices and
$\tr(1_{n})=n$; in the second case, $\lambda= 0$, we assume that
$\tr(1)=1$.

Observe that, for any integer $\lambda=n$, we have
$$
\frac{e^{nt}-e^{-nt}}{e^{t}-e^{-t}}=e^{(1-n)t}+e^{(3-n)t}+\dots
+e^{(n-3)t}+e^{(n-1)t}
$$
and the respective functional is of the form
$$
\tr(f(H))=\mathop{\sum}\limits_{i=1}^nf(\alpha_{i}), \quad
\alpha_{i}=n-2i+1.
$$

The functional $\tr$ on $\fgl(\lambda)$ gives rise to a
non-degenerate symmetric invariant bilinear form $\langle u,
v\rangle=\tr~uv$.

\ssbegin{9.1}{Theorem} In notations of \S $8$: {\em i)} $\langle e_{k,
l}, e_{k', l'}\rangle=\delta_{k, k'}\delta_{l+l', 0}$.

{\em ii)} The polynomials $f_{k, l}$ for a fixed $l\geq 0$ and $k\geq l$
form an orthogonal basis with respect to the scalar product
$(8.1.1)$.

{\em iii)} The polynomials $f_{k, l}$ coincide, up to a constant
factor, with continuous Hahn polynomials of discrete variable
given by $(8.1.2)$ for $n=\lambda$.\end{Theorem}

\begin{proof} The proof of headings i) and ii) is similar to that of
respective statements of Theorem 8.1 for integer $\lambda$.  To
prove iii), observe that, for a nonnegative integer, eq.  (9.1)
holds thanks to Theorem 8.1.  But both sides of (9.1) are
polynomials in $H$ whose coefficients are polynomials in
$\lambda$.  Due to continuity of these expressions in Zariski
topology we are done.
\end{proof}

\ssbegin{9.2}{Theorem} Let $\Omega$ be the quadratic Casimir operator
for $\fsl(2)\subset\fgl(\lambda)$.

{\em i)} $\Omega$ is self-adjoint with respect to the form $(8.1.1)$
and the polynomials $X^lf_{kl}$ for $k\in\Zee_{+}$ and $l\in[-k,
k]\cap\Zee$ are eigenfunctions of $\Omega$ corresponding to eigenvalue
$2k(k+1)$.  The polynomials $f_{kl}$ satisfy the difference equation
$(8.2.i)$ and are of the form $(8.2.ii)$

{\em ii)} $\langle f_{kl}, f_{kl}\rangle=
\begin{cases}\displaystyle\frac{(k-l)!}{(k+l)!}
\displaystyle\frac{(k!)^2}{2k+1}\lambda(\lambda^2-1^{2})\dots
(\lambda^2-k^{2})&\text{if $\lambda\neq 0$}\\
(-1)^k\displaystyle\frac{(k-l)!}{(k+l)!}\displaystyle\frac{(k!)^4}{2k+1}&
\text{if $\lambda=0$.} \end{cases}$
\end{Theorem}

\begin{proof} Proof of heading i) is the same as that of the
corresponding statements of Theorem 8.2.  If $\lambda\neq 0$, then
both sides of eq.  ii) are polynomials in $\lambda$ which coincide
at integer values of $\lambda$ and, therefore, are equal by
continuity in Zariski topology.  To embrace $\lambda=0$, consider
$\mathop{\lim}\limits_{\lambda\tto
0}\displaystyle\frac{\tr}{\lambda}$; as before, both sides become
polynomials in $\lambda$ implying ii).  \end{proof}

It is somewhat unexpected that the dual orthogonality relations hold
for sufficiently large $|\lambda|$ .  Indeed, set
$$
c_{kl}=\begin{cases}\displaystyle\frac{(k-l)!}{(k+l)!}
\displaystyle\frac{(k!)^2}{2k+1}(\lambda^2-1^{2})\dots
(\lambda^2-k^{2})&\text{ if $\lambda\neq 0$}\\
(-1)^k\displaystyle\frac{(k-l)!}{(k+l)!}\displaystyle\frac{(k!)^4}{2k+1}&
\text{ if $\lambda=0$.} \end{cases}\eqno{(9.2)}
$$

\ssbegin{9.3}{Conjecture} For sufficiently large $|\lambda|$, we
have (for $\alpha_{i}=\lambda-2i+1$):
$$
\mathop{\sum}\limits_{k=l}^{\infty}\frac{1}{c_{kl}}f_{kl}(\alpha_{i})
f_{kl}(\alpha_{j}) T_{1}(\alpha_{i})\dots
T_{l}(\alpha_{i})=\delta_{ij} \text{ for $i, j>l$}.\eqno{(9.3.1)}
$$
$$
\mathop{\sum}\limits_{k=0}^{\infty}\left(
\frac{1}{c_{k0}}f_{k0}(\alpha_{i})^2+2\mathop{\sum}
\limits_{l=1}^{k}\frac{1}{c_{kl}} f_{kl}(\alpha_{i})^2
T_{1}(\alpha_{i})\dots
T_{l}(\alpha_{i})-\displaystyle\frac{2k+1}{\lambda}\right)=
-\alpha_{i}\text{ for $i\in\Nee$}.\eqno{(9.3.2)}
$$
\end{Conjecture}

A priori, $|\lambda|$ depends on $l$. We do not know a uniform proof
and only checked Conjecture in certain particular cases. We are
thankful to V.~Gerdt and P.~Grozman who helped us with numerical
experiments on Maple and {\it Mathematica}, respectively.

\ssec{9.4.  On positive definiteness of the form $\langle \cdot ,
\cdot \rangle$} Formulas above make it clear that if we divide the
form $\langle \cdot , \cdot \rangle$ by $\lambda$ and consider
polynomials of even degree only we get a sign-definite form not
only for integer values of $\lambda$ but also for real values such
that $0<|\lambda|<1$ and for purely imaginary $\lambda$.  The last
observation suggests to divide the trace by $\lambda$ from the
very beginning.

\end{document}